\documentclass{amsart}

\usepackage[leqno]{amsmath}
\usepackage{amsthm}
\usepackage{latexsym,amsfonts,amssymb}
\usepackage[all]{xy} \SelectTips{eu}{}
\usepackage{hyperref}

\newcommand{\numberseries}{\bfseries}   
\newlength{\thmtopspace}                
\newlength{\thmbotspace}                
\newlength{\thmheadspace}               
\newlength{\thmindent}                  
\newtheoremstyle{bfupright head,slanted body}
                {\thmtopspace}{\thmbotspace}
                {\slshape}{\thmindent}{\bfseries}{.}{\thmheadspace}
                {{\numberseries \thmnumber{#2\;}}\thmnote{#3}}
\newtheoremstyle{fixed bf head,slanted body}
                {\thmtopspace}{\thmbotspace}{\slshape}
                {\thmindent}{\bfseries}{.}{\thmheadspace}
                {{\numberseries \thmnumber{#2\;}}\thmname{#1}\thmnote{ (#3)}}

\newtheoremstyle{fixed bf head,upright body}
                {\thmtopspace}{\thmbotspace}{\upshape}
                {\thmindent}{\bfseries}{.}{\thmheadspace}
                {{\numberseries \thmnumber{#2\;}}\thmname{#1}\thmnote{ (#3)}}
\newtheoremstyle{numbered paragraph}
                {\thmtopspace}{\thmbotspace}{\upshape}
                {\thmindent}{\upshape}{}{\thmheadspace}
                {{\numberseries \thmnumber{#2.}}}
\theoremstyle{bfupright head,slanted body}
\newtheorem{res}{}[section]             \newtheorem*{res*}{}
\theoremstyle{fixed bf head,slanted body}
\newtheorem{thm}[res]{Theorem}          \newtheorem*{thm*}{Theorem}
\newtheorem{prp}[res]{Proposition}      \newtheorem*{prp*}{Proposition}
\newtheorem{cor}[res]{Corollary}        \newtheorem*{cor*}{Corollary}
\newtheorem{lem}[res]{Lemma}            \newtheorem*{lem*}{Lemma}
\theoremstyle{fixed bf head,upright body}
\newtheorem{dfn}[res]{Definition}       \newtheorem*{dfn*}{Definition}
\newtheorem{rmk}[res]{Remark}           \newtheorem*{rmk*}{Remark}
          \newtheorem*{exa*}{Example}
\theoremstyle{numbered paragraph}
\newtheorem{ipg}[res]{}
\newlength{\thmlistleft}        
\newlength{\thmlistright}       
\newlength{\thmlistpartopsep}   
\newlength{\thmlisttopsep}      
\newlength{\thmlistparsep}      
\newlength{\thmlistitemsep}     
\newcounter{eqc}
\newenvironment{eqc}{\begin{list}{\upshape (\textit{\roman{eqc}})}%
    {\usecounter{eqc}%
      \setlength{\leftmargin}{\thmlistleft}%
      \setlength{\labelwidth}{\thmlistleft}%
      \setlength{\rightmargin}{\thmlistright}%
      \setlength{\partopsep}{\thmlistpartopsep}%
      \setlength{\topsep}{\thmlisttopsep}%
      \setlength{\parsep}{\thmlistparsep}%
      \setlength{\itemsep}{\thmlistitemsep}}}%
  {\end{list}}%
\newcommand{\eqclbl}[1]{{\upshape(\textit{#1})}}
\newcounter{prt}
\newenvironment{prt}{\begin{list}{\upshape (\alph{prt})}%
    {\usecounter{prt}%
      \setlength{\leftmargin}{\thmlistleft}%
      \setlength{\labelwidth}{\thmlistleft}%
      \setlength{\rightmargin}{\thmlistright}%
      \setlength{\partopsep}{\thmlistpartopsep}%
      \setlength{\topsep}{\thmlisttopsep}%
      \setlength{\parsep}{\thmlistparsep}%
      \setlength{\itemsep}{\thmlistitemsep}}}%
  {\end{list}}%
\newcommand{\prtlbl}[1]{{\upshape(#1)}}
\newcounter{rqm}
\newenvironment{rqm}{\begin{list}{\upshape (\arabic{rqm})}%
    {\usecounter{rqm}%
      \setlength{\leftmargin}{\thmlistleft}%
      \setlength{\labelwidth}{\thmlistleft}%
      \setlength{\rightmargin}{\thmlistright}%
      \setlength{\partopsep}{\thmlistpartopsep}%
      \setlength{\topsep}{\thmlisttopsep}%
      \setlength{\parsep}{\thmlistparsep}%
      \setlength{\itemsep}{\thmlistitemsep}}}%
  {\end{list}}%
\newcommand{\rqmlbl}[1]{{\upshape(#1)}}
\newenvironment{itemlist}{\nopagebreak \begin{list}{$\bullet$}%
    {\setlength{\leftmargin}{1.5em}%
      \setlength{\labelwidth}{\thmlistleft}%
      \setlength{\rightmargin}{\thmlistright}%
      \setlength{\partopsep}{\thmlistpartopsep}%
      \setlength{\topsep}{\thmlisttopsep}%
      \setlength{\parsep}{\thmlistparsep}%
      \setlength{\itemsep}{\thmlistitemsep}}}%
  {\end{list}}%
\newenvironment{prf*}[1][Proof]{%
  \begin{proof}[\bf #1]
    \setcounter{equation}{0}
    }
  {\end{proof}
}
\newcommand{\proofofimp}[3][:]{\mbox{\eqclbl{#2}$\!\implies\!$\eqclbl{#3}#1}}

\newcommand{\pgref}[1]{\ref{#1}}
\newcommand{\thmref}[2][Theorem~]{#1\pgref{thm:#2}}
\newcommand{\corref}[2][Corollary~]{#1\pgref{cor:#2}}
\newcommand{\prpref}[2][Proposition~]{#1\pgref{prp:#2}}
\newcommand{\lemref}[2][Lemma~]{#1\pgref{lem:#2}}
\newcommand{\dfnref}[2][Definition~]{#1\pgref{dfn:#2}}
\newcommand{\rmkref}[2][Remark~]{#1\pgref{rmk:#2}}
\newcommand{\secref}[2][Section~]{#1\ref{sec:#2}}
\newcommand{\partpgref}[2]{\ref{#1}\prtlbl{#2}}
\newcommand{\partprpref}[3][Proposition~]{#1\partpgref{prp:#2}{#3}}
\renewcommand{\eqref}[1]{(\pgref{eq:#1})}
\newcommand{\thmcite}[2][?]{\cite[thm.~#1]{#2}}
\newcommand{\corcite}[2][?]{\cite[cor.~#1]{#2}}
\newcommand{\prpcite}[2][?]{\cite[prop.~#1]{#2}}
\newcommand{\lemcite}[2][?]{\cite[lem.~#1]{#2}}

\newcommand{\seccite}[2][?]{\cite[sec.~#1]{#2}}
\numberwithin{equation}{res}

\newcommand{\fsp}[1][1.2]{f\hspace{#1pt}}
\newcommand{\dimR}{\operatorname{dim}R}
\newcommand{\setof}[3][\mspace{1mu}]{\{#1#2 \mid #3#1\}}
\newcommand{\kk}{\Bbbk}
\newcommand{\NN}{\mathbb{N}}
\newcommand{\ZZ}{\mathbb{Z}}
\newcommand{\qtext}[1]{\quad\text{#1}\quad}
\newcommand{\qqtext}[1]{\qquad\text{#1}\qquad}
\newcommand{\qand}{\qtext{and}}
\newcommand{\qqand}{\qqtext{and}}
\newcommand{\deq}{\:=\:}
\newcommand{\dis}{\:\is\:}
\newcommand{\dqis}{\:\qis\:}
\renewcommand{\a}{\alpha}
\newcommand{\m}{\mathfrak{m}}
\newcommand{\p}{\mathfrak{p}}

\newcommand{\is}{\cong}
\newcommand{\qis}{\simeq}
\renewcommand{\le}{\leqslant}
\renewcommand{\ge}{\geqslant}
\newcommand{\lra}{\longrightarrow}
\newcommand{\xra}[2][]{\xrightarrow[#1]{\:#2\:}}
\newcommand{\qra}{\xra{\smash{\qis}}}
\newcommand{\Rop}{R^\circ}
\newcommand{\Sop}{S^\circ}
\newcommand{\mapdef}[4][\rightarrow]{\nobreak{#2\colon #3 #1 #4}}
\newcommand{\dmapdef}[4][\lra]{\nobreak{#2\colon #3\:#1\:#4}}
\newcommand{\Ker}[1]{\nobreak{\operatorname{Ker}#1}}
\newcommand{\Coker}[1]{\nobreak{\operatorname{Coker}#1}}
\renewcommand{\Im}[1]{\nobreak{\operatorname{Im}#1}}
\newcommand{\Cone}[1]{\nobreak{\operatorname{Cone}#1}}
\newcommand{\dif}[2][]{{\partial}^{#2}_{#1}}
\newcommand{\Co}[2][]{\operatorname{C}_{#1}(#2)}
\renewcommand{\H}[2][]{\operatorname{H}_{#1}(#2)}
\newcommand{\Tha}[2]{#2_{{\scriptscriptstyle\le}#1}}
\newcommand{\Thb}[2]{#2_{{\scriptscriptstyle\ge}#1}}
\newcommand{\Tsb}[2]{#2_{{\scriptscriptstyle\supset}#1}}

\newcommand{\dpt}[2][R]{\operatorname{depth}_{#1}#2}
\newcommand{\fd}[2][R]{\operatorname{fd}_{#1}#2}
\newcommand{\id}[2][R]{\operatorname{id}_{#1}#2}
\newcommand{\pd}[2][R]{\operatorname{pd}_{#1}#2}
\newcommand{\Gdim}[2][R]{\operatorname{G-dim}_{#1}#2}
\newcommand{\Hom}[3][R]{\operatorname{Hom}_{#1}(#2,#3)}
\newcommand{\RHom}[3][R]{\operatorname{\mathbf{R}Hom}_{#1}(#2,#3)}
\newcommand{\Ext}[4][R]{\operatorname{Ext}_{#1}^{#2}(#3,#4)}
\newcommand{\tp}[3][R]{\nobreak{#2\otimes_{#1}#3}}
\newcommand{\tpp}[3][R]{(\tp[#1]{#2}{#3})}

\newcommand{\Tor}[4][R]{\operatorname{Tor}^{#1}_{#2}(#3,#4)}
\newcommand{\Cat}[2]{{\mathsf{#2}}(#1)}
\newcommand{\Mod}[1][R]{\Cat{#1}{M}}
\newcommand{\D}[1][R]{\Cat{#1}{D}}
\newcommand{\cone}[2]{\operatorname{C}^{#1}_{#2}}
\newcommand{\Dual}[1]{\operatorname{D}(#1)}
\newcommand{\dcx}[3]{#1_{#2,#3}}
\newcommand{\del}{\partial}

\newcommand{\btp}[3][R]{\nobreak{#2\mathbin{\widebar{\otimes}}_{#1}#3}}
\newcommand{\btpp}[3][R]{(\nobreak{#2\mathbin{\widebar{\otimes}}_{#1}#3})}
\newcommand{\ttp}[3][R]{\nobreak{#2\mathbin{\widetilde{\otimes}}_{#1}#3}}
\newcommand{\ttpp}[3][R]{(\nobreak{#2\mathbin{\widetilde{\otimes}}_{#1}#3})}
\newcommand{\ptp}[3][R]{\nobreak{#2\otimes^\Join_{#1}#3}}
\newcommand{\ptpp}[3][R]{\nobreak{(#2\otimes^\Join_{#1}#3)}}
\newcommand{\bTor}[4][R]{\widebar{\operatorname{Tor}}_{#2}^{#1}(#3,#4)}
\newcommand{\Ttor}[4][R]{\smash{\operatorname{\widehat{Tor}}}_%
  {#2}^{{#1}^{\phantom{|\mspace{-6mu}}}}(#3,#4)}
\newcommand{\Stor}[4][R]{\smash{\operatorname{\widetilde{Tor}}}_%
  {#2}^{{#1}^{\phantom{|\mspace{-6mu}}}}(#3,#4)}

\newcommand{\bExt}[4][R]{\widebar{\operatorname{Ext}}_{#1}^{#2}(#3,#4)}
\newcommand{\Text}[4][R]{\smash{\operatorname{\widehat{Ext}}%
  }_{#1}^{#2^{\phantom{|}\mspace{-6mu}}}(#3,#4)}
\newcommand{\Sext}[4][R]{\smash{\operatorname{\widetilde{Ext}}%
  }_{#1}^{#2^{\phantom{|}\mspace{-6mu}}}(#3,#4)}
\newcommand{\bHom}[3][R]{\widebar{\operatorname{Hom}}_{#1}(#2,#3)}
\newcommand{\SHom}[3][R]{\widetilde{\operatorname{Hom}}_{#1}(#2,#3)}
\newcommand{\susp}{\operatorname{\mathsf{\Sigma}}}
\newcommand{\Hh}[2][]{\operatorname{H}^h_{#1}(#2)}
\newcommand{\Hv}[2][]{\operatorname{H}^v_{#1}(#2)}
\newcommand{\rank}{\operatorname{rank}}
\DeclareMathOperator*{\dcoprod}{\textstyle\coprod}
\DeclareMathOperator*{\dprod}{\textstyle\prod}
\def\urltilda{\kern
  -.15em\lower .7ex\hbox{\~{}}\kern .04em} \makeatletter
\def\@nobreak@#1{\mathchoice%
  {\nobreakdef@\displaystyle\f@size{#1}}%
  {\nobreakdef@\nobreakstyle\tf@size{\firstchoice@false #1}}%
  {\nobreakdef@\nobreakstyle\sf@size{\firstchoice@false #1}}%
  {\nobreakdef@\nobreakstyle\ssf@size{\firstchoice@false #1}}%
  \check@mathfonts}%
\def\nobreakdef@#1#2#3{\hbox{{%
      \everymath{#1}%
      \let\f@size#2\selectfont%
      #3}}}%
\makeatother

\def\widebardisplay#1{%
  \setbox0=\hbox{$\displaystyle #1$} \dimen0=\wd0%
  \advance\dimen0 by -1pt
  \vbox{%
    \nointerlineskip%
    \moveright
    .25pt 
    \vbox{\hrule width \dimen0}%
    \nointerlineskip%
    \kern 2.25pt
    \box0%
  }%
}

\def\widebartext#1{%
  \setbox0=\hbox{$#1$} \dimen0=\wd0%
  \advance\dimen0 by -1pt
  \vbox{%
    \nointerlineskip%
    \moveright
    .25pt 
    \vbox{\hrule width \dimen0}%
    \nointerlineskip%
    \kern 1.75pt
    \box0%
  }%
}

\def\widebarscript#1{%
  \setbox0=\hbox{$\scriptstyle #1$} \dimen0=\wd0%
  \advance\dimen0 by -2pt
  \vbox{%
    \nointerlineskip%
    \moveright
    1pt 
    \vbox{\hrule width \dimen0}%
    \nointerlineskip%
    \kern .8pt
    \box0%
  }%
}

\def\widebarscriptscript#1{%
  \setbox0=\hbox{$\scriptscriptstyle #1$} \dimen0=\wd0%
  \advance\dimen0 by -2pt
  \vbox{%
    \nointerlineskip%
    \moveright
    1pt 
    \vbox{\hrule width \dimen0}%
    \nointerlineskip%
    \kern .6pt
    \box0%
  }%
}

\def\widebar#1{\mathchoice%
  {\widebardisplay{#1}}%
  {\widebartext{#1}}%
  {\widebarscript{#1}}%
  {\widebarscriptscript{#1}}%
}

\DeclareFontFamily{U}{mathx}{\hyphenchar\font45}
\DeclareFontShape{U}{mathx}{m}{n}{ <5> <6> <7> <8> <9> <10> <10.95>
  <12> <14.4> <17.28> <20.74> <24.88> mathx10 }{}
\DeclareSymbolFont{mathx}{U}{mathx}{m}{n}
\DeclareFontSubstitution{U}{mathx}{m}{n}
\DeclareMathAccent{\widecheck}{0}{mathx}{"71}
\DeclareMathAccent{\wideparen}{0}{mathx}{"75}

\begin{document}

\title{Stable homology over associative rings}

\author[O. Celikbas]{Olgur Celikbas}

\address{University of Missouri, Columbia, MO~65211, U.S.A.}

\curraddr{University of Connecticut, Storrs, CT~06269, U.S.A}

\email{olgur.celikbas@uconn.edu}

\author[L.\,W. Christensen]{Lars Winther Christensen}

\address{Texas Tech University, Lubbock, TX 79409, U.S.A.}

\email{lars.w.christensen@ttu.edu}

\urladdr{http://www.math.ttu.edu/\urltilda lchriste}

\author[L. Liang]{Li Liang}

\address{Lanzhou Jiaotong University, Lanzhou 730070, China}

\email{lliangnju@gmail.com}

\author[G. Piepmeyer]{Greg Piepmeyer}

\address{Columbia Basin College, Pasco, WA~99301, U.S.A.}

\email{pggreg@gmail.com}

\thanks{This research was partly supported by a Simons Foundation
  Collaboration Grant for Mathematicians, award no. 281886 (L.W.C.),
  NSA grant H98230-14-0140 (L.W.C.) and NSFC grant 11301240
  (L.L.). Part of the work was done during L.L.'s stay at Texas Tech
  University with support from the China Scholarship Council; he
  thanks the Department of Mathematics and Statistics at Texas Tech
  for its kind hospitality.}

\date{29 December 2015}





\begin{abstract}
  We analyze \emph{stable homology} over associative rings and obtain
  results over Artin algebras and commutative noetherian rings.  Our
  study develops similarly for these classes; for simplicity we only
  discuss the latter here.

  Stable homology is a broad generalization of Tate
  homology. Vanishing of stable homology detects classes of
  rings---among them Gorenstein rings, the original domain of Tate
  homology. Closely related to gorensteinness of rings is
  Auslander's G-dimension for modules. We show that vanishing of
  stable homology detects modules of finite G-dimension. This is the
  first characterization of such modules in terms of vanishing of
  (co)homology \textsl{alone.}

  Stable homology, like absolute homology, Tor, is a theory in two
  variables. It can be computed from a flat resolution of one module
  together with an injective resolution of the other. This betrays
  that stable homology is not balanced in the way Tor is balanced. In
  fact, we prove that a ring is Gorenstein if and only if stable homology is
  balanced.
\end{abstract}

\maketitle

\thispagestyle{empty}

\section*{Introduction}

\noindent
The homology theory studied in this paper was introduced by P.~Vogel
in the 1980s. Vogel did not publish his work, but the theory appeared
in print in a 1992 paper by Goichot \cite{FGc92}, who called it
\emph{Tate--Vogel homology.} As the name suggests, the theory is a
generalization of Tate homology for modules over finite group
rings. Vogel and Goichot also considered a generalization of Tate
cohomology, which was studied in detail by Avramov and
Veliche~\cite{LLAOVl07}. In that paper, the theory was called
\emph{stable cohomology}, to emphasize a relation to stabilization of
module categories. To align terminology, we henceforth refer to the
homology theory as \emph{stable homology}.

For modules $M$ and $N$ over a ring $R$, stable homology is a
$\ZZ$-indexed family of abelian groups $\Stor{i}{M}{N}$. These fit
into an exact sequence
\begin{equation*}
  \label{dagger}
  \tag{$\dag$}
  \cdots \to \Stor{i}{M}{N} \to \Tor{i}{M}{N} \to \bTor{i}{M}{N} \to
  \Stor{i-1}{M}{N} \to \cdots
\end{equation*}
where the groups $\bTor{i}{M}{N}$ and $\Tor{i}{M}{N}$ form,
respectively, the \emph{unbounded homology} and the standard
\emph{absolute homology} of $M$ and $N$. We are thus led to study
stable and unbounded homology simultaneously. Our investigation takes
cues from the studies of stable cohomology \cite{LLAOVl07} and
absolute homology, but our results look quite different.  This comes
down to an inherent asymmetry in the definition of stable homology
that is not present in either of these precursors. It manifests itself
in different ways, but it is apparent in most of our results.
\begin{equation*}
  \ast \ \ \ast \ \ \ast
\end{equation*}
When we consider stable homology $\Stor{}{M}{N}$, the ring acts on $M$
from the right\linebreak and on $N$ from the left. In this paper an
$R$-module is a left $R$-module; to distinguish\linebreak right
$R$-modules, we speak of modules over the opposite ring $\Rop$.  We
study stable homology over associative rings and obtain conclusive
results for Artin algebras and commutative noetherian rings. In the
following overview, $R$ denotes an Artin algebra or a commutative
noetherian local ring.  The results we discuss are special cases of
results obtained within the paper; internal references are given in
parentheses.

One expects a homology theory to detect finiteness of homological
dimensions. Our first two results reflect the asymmetry in the
definition of stable homology.

\begin{res*}[Right vanishing {\rm(\thmref[]{finitepd},
    \corref[]{artin pd})}]
  For a finitely generated $\Rop$-module $M$, the following conditions
  are equivalent.
  \begin{eqc}
  \item $M$ has finite projective dimension.
  \item $\Stor{i}{M}{-}=0$ for all $i\in\ZZ$.
  \item There is an $i \ge 0$ with $\Stor{i}{M}{-}=0$.
  \end{eqc}
\end{res*}

\begin{res*}[Left vanishing {\rm(\prpref[]{artin id},
    \corref[]{finiteid})}]
  \label{leftvanishing}
  For a finitely generated $R$-module $N$, the following conditions
  are equivalent.
  \begin{eqc}
  \item $N$ has finite injective dimension.
  \item $\Stor{i}{-}{N}=0$ for all $i\in\ZZ$.
  \item There is an $i \le 0$ with $\Stor{i}{-}{N}=0$.
  \end{eqc}
\end{res*}

\noindent
These two vanishing results reveal that stable homology cannot be
balanced in the way absolute homology, Tor, is balanced.

\begin{res*}[Balancedness {\rm(\corref[]{art-balance},
    \corref[]{Gor2}, \corref[]{Gor3})}] The following conditions on $R$ are
  equivalent.
  \begin{eqc}
  \item $R$ has finite injective dimension over $R$ and over $\Rop$.
  \item For all finitely generated $\Rop$-modules $M$, all finitely
    generated $R$-modules $N$, and all $i\in\ZZ$ there are
    isomorphisms $\Stor{i}{M}{N} \is \Stor[\Rop]{i}{N}{M}$.
  \item For all $\Rop$-modules $M$, all $R$-modules $N$, and all
    $i\in\ZZ$ there are isomorphisms $\Stor{i}{M}{N} \is
    \Stor[\Rop]{i}{N}{M}$.
  \end{eqc}
\end{res*}

\noindent
Another way to phrase part \eqclbl{i} above is to say that $R$ is
\emph{Iwanaga-Gorenstein}. On that topic: a commutative noetherian
ring $A$ is \emph{Gorenstein} (\emph{regular} or
\emph{Cohen--Macaulay}) if the local ring $A_\p$ is Iwanaga-Gorenstein
(regular or Cohen--Macaulay) for every prime ideal $\p$ in $A$. It
follows that a commutative noetherian ring of finite self-injective
dimension is Gorenstein, but a Gorenstein ring need not have finite
self-injective dimension; consider, for example, Nagata's regular ring
of infinite Krull dimension \cite[appn.\ exa.~1]{Nag}.

It came as a surprise to us that balancedness of stable homology
detects gorensteinness outside of the local situation
(\corref[]{Gor2}): A commutative noetherian ring $A$ is Gorenstein if
and only if $\Stor[A]{i}{M}{N}$ and $\Stor[A]{i}{N}{M}$ are isomorphic
for all finitely generated $A$-modules $M$ and $N$ and all
$i\in\ZZ$. This turns out to be only one of several ways in which
stable homology captures global properties of rings. Indeed, vanishing
of stable homology detects if a commutative noetherian (not
necessarily local) ring is regular~(\corref[]{regular}), Gorenstein
(\corref[]{Gor1}), or Cohen--Macaulay (\corref[]{cm}).

For finitely generated modules over a noetherian ring, Auslander and
Bridger~\cite{MAsMBr69} introduced a homological invariant called the
\emph{G-dimension}. A noetherian ring is Iwanaga-Gorenstein if and
only if there is an integer $d$ such that every finitely generated
left or right module has G-dimension at most $d$; see Avramov and
Martsinkovsky~\cite[3.2]{LLAAMr02}. In \seccite[3]{FGc92} Tate
homology was defined for modules over Iwanaga-Gorenstein rings and
shown to agree with stable homology. Iacob~\cite{AIc07} generalized
Tate homology further to a setting that includes the case where $M$ is
a finitely generated $\Rop$-module of finite G-dimension. We prove:

\begin{res*}[Tate homology {\rm (\thmref[]{Gproj-Gflat-fg})}]
  Let $M$ be a finitely generated $\Rop$-module of finite
  G-dimension. For every $i\in\ZZ$ the stable homology functor
  $\Stor{i}{M}{-}$ is isomorphic to the Tate homology functor
  $\Ttor{i}{M}{-}$.
\end{res*}

\noindent
There are two primary generalizations of G-dimension to modules that
are not finitely generated: \emph{Gorenstein projective dimension} and
\emph{Gorenstein flat dimension.} They are defined in terms of
resolutions by Gorenstein projective (flat) modules, notions which were
introduced by Enochs, Jenda, and collaborators; see Holm
\cite{HHl04a}\footnote{ Enochs and Jenda consolidated their work in
  \cite[chaps.~10--12]{rha}, which deal almost exclusively with
  Gorenstein rings. In \cite{HHl04a} Holm does Gorenstein homological
  algebra over associative rings, so that is our standard reference
  for this topic. Basic definitions are recalled in \pgref{ipg:gpgf}.}.
In general, it is not known if Gorenstein projective modules are
Gorenstein flat.

Iacob's definition of Tate homology $\Ttor{}{M}{-}$ is for
$\Rop$-modules $M$ of finite Gorenstein projective dimension. We show
(\thmref[]{Gproj-Gflat}) that Tate homology agrees with stable
homology $\Stor{}{M}{-}$ for all such $\Rop$-modules if and only if
every Gorenstein projective $\Rop$-module is Gorenstein flat.

Here is an example to illustrate how widely the various homology
theories differ.

\begin{exa*}
  Let $(R,\m)$ be a commutative artinian local ring with $\m^2 =0$, and
  assume that $R$ is not Gorenstein. If $k$ denotes a field, then
  $k[x,y]/(x^2, xy, y^2)$ is an example of such a ring. Let $E$ be the
  injective hull of the residue field $R/\m$. Then: 
  \begin{prt}
  \item Absolute homology $\Tor{i}{E}{E}$ is non-zero for every $i\geq
    2$.

  \item Stable homology $\Stor{i}{E}{E}$ is zero for every $i$.

  \item Unstable homology $\bTor{i}{E}{E}$ agrees with $\Tor{i}{E}{E}$
    for every $i$.

  \item Tate homology ``$\Ttor{}{E}{E}$'' is not defined.
  \end{prt}
  Since $R$ is not Gorenstein, $E$ has infinite G-dimension; see
  \thmcite[3.2]{HHl04c}. That explains (d); see \ref{ipg:gproj} and
  \pgref{ipg:gp}. It also explains (a), as the first syzygy of $E$ is
  a $k$-vector space, whence vanishing of $\Tor{i}{E}{E}$ for some $i
  \ge 2$ would imply that $E$ has finite projective dimension and
  hence finite G-dimension; see \thmcite[4.9]{LLAAMr02}. Left
  vanishing, see p.\:\pageref{leftvanishing}, yields (b), and then the
  exact sequence $(\dag)$ on p.\:\pageref{dagger} gives~(c).
\end{exa*}

Finally we return to the topic of vanishing. A much studied question
in Gorenstein homological algebra is how to detect finiteness of
G-dimension. A finitely generated module $G$ over a commutative
noetherian ring $A$ has G-dimension zero if it satisfies
$\Ext[A]{i}{G}{A} = 0 = \Ext[A]{i}{\Hom[A]{G}{A}}{A}$ for all $i>0$
\textsl{and} the canonical map $G \to \Hom[A]{\Hom[A]{G}{A}}{A}$ is an
isomorphism.  This is the original definition~\cite{MAsMBr69}, and the
last requirement cannot be dispensed with, as shown by Jorgensen and
\c{S}ega~\cite{DAJLMS06}. There are other characterizations of modules
of finite G-dimension in the literature, but none that can be
expressed solely in terms of vanishing of (co)homology functors.
Hence, we were excited to discover:

\begin{res*}
[G-dimension {\rm (\thmref[]{vanishing})}]
Let $A$ be a commutative noetherian ring. A finitely generated
$A$-module $M$ has finite G-dimension if and only if\,
\mbox{$\Stor[A]{i}{M}{A} = 0$} holds for all $i\in\ZZ$.
\end{res*}

\section{Tensor products of complexes}
\label{sec:tensor}

\noindent
All rings in this paper are assumed to be associative algebras over a
commutative ring $\kk$, where $\kk = \ZZ$ is always possible, and $\kk
= R$ works when $R$ is commutative. For an Artin algebra $R$, one can
take $\kk$ to be an artinian ring such that $R$ is finitely generated
as a $\kk$-module. We recall that in this situation the functor
$\Dual{-} = \Hom[\kk]{-}{E}$, where $E$ is the injective hull of
$\kk/\mathrm{Jac}(\kk)$, yields a duality between the categories of
finitely generated $R$-modules and finitely generated $\Rop$-modules.

By $\Mod$ we denote the category of $R$-modules. We work with
complexes, index these homologically, and follow standard notation
(see the appendix in \cite{lnm} for what is not covered below).

\begin{ipg}
  \label{intro notation}
  For a complex $X$ with differential $\dif[]{X}$ and an integer
  $n\in\ZZ$ the symbol $\Co[n]{X}$ denotes the cokernel of
  $\dif[n+1]{X}$, and $\H[n]{X}$ denotes the homology of $X$ in degree
  $n$, i.e., $\Ker{\dif[n]{X}}/\Im{\dif[n+1]{X}}$. A complex with
  $\H{X}=0$ is called \emph{acyclic,} and a morphism of complexes $X
  \to Y$ that induces an isomorphism $\H{X} \to \H{Y}$ is called a
  \emph{quasi-isomorphism}. The symbol $\qis$ is used to decorate
  quasi-isomorphisms; it is also used for isomorphisms in derived
  categories.
\end{ipg}

\begin{ipg}
  \label{std hom and tensor}
  For an $\Rop$-complex $X$ and an $R$-complex $Y$, the tensor product
  $\tp{X}{Y}$ is the $\kk$-complex with degree $n$ term
  $(X\otimes_{R}Y)_{n}=\coprod_{i\in\ZZ}(X_{i}\otimes_{R}Y_{n-i})$ and
  differential given by $\del(x\otimes y)=\dif[i]{X}(x)\otimes
  y+(-1)^{i}x\otimes\dif[n-i]{Y}(y)$ for $x\in X_{i}$ and $y\in
  Y_{n-i}$.
\end{ipg}

Contrast the standard tensor product complex in \ref{std hom and
  tensor} with the construction in \ref{stable and unbounded tensor
  product}, which first appeared in \cite{FGc92}. Similar
constructions for Hom were also treated in \cite{FGc92} and in great
detail by Avramov and Veliche~\cite{LLAOVl07}. We recall the Hom
constructions in the Appendix, where we also study their interactions
with the one below.

\begin{dfn}
  \label{stable and unbounded tensor product}
  For an $\Rop$-complex $X$ and an $R$-complex $Y$, consider the
  graded $\kk$-module $\btp{X}{Y}$ with degree $n$ component
  \begin{equation*}
    (\btp{X}{Y})_n =
    \prod_{i\in\ZZ}\tp{X_i}{Y_{n-i}}\:.
  \end{equation*}
  Endowed with the degree $-1$ homomorphism $\del$ defined on
  elementary tensors as in \ref{std hom and tensor}, it becomes a
  complex called the \emph{unbounded} tensor product. It contains the
  tensor product, $\tp{X}{Y}$, as a subcomplex. The quotient complex
  $(\btp{X}{Y})/(\tp{X}{Y})$ is denoted $\ttp{X}{Y}$, and it is called
  the \emph{stable} tensor product.
\end{dfn}

\begin{ipg}
  \label{simplifications}
  The definition above yields an exact sequence of $\kk$-complexes,
  \begin{equation}
    \label{eq:ts}
    0 \lra \tp{X}{Y} \lra \btp{X}{Y} \lra \ttp{X}{Y} \lra 0\:.
  \end{equation}
  Notice that if $X$ or $Y$ is bounded, or if both complexes are
  bounded on the same side (above or below), then the unbounded tensor
  product coincides with the usual tensor product, and the stable
  tensor product $\ttp{X}{Y}$ is zero.

  Note that the unbounded tensor product $\btp{X}{Y}$ is the
  product totalization of the double complex $\tpp[R]{X_i}{Y_j}$.
\end{ipg}

\begin{ipg}
  \label{btp}
  We collect some basic results about the unbounded and stable tensor
  products. For the unbounded tensor
  product the proofs mimic the proofs for the tensor product, and one
  obtains the results for the stable tensor product via \eqref{ts}. As
  above $X$ is an $\Rop$-complex and $Y$ is an $R$-complex.
  \begin{prt}
  \item There are isomorphisms of $\kk$-complexes
    \begin{equation*}
      \btp{X}{Y} \dis \btp[\Rop]{Y}{X} \qqand \ttp{X}{Y} \dis
      \ttp[\Rop]{Y}{X}\:.
    \end{equation*}
  \item The functors $\btp{X}{-}$ and $\ttp{X}{-}$ are additive and
    right exact.
  \item The functors $\btp{X}{-}$ and $\ttp{X}{-}$ preserve
    degree-wise split exact sequences.
  \item The functors $\btp{X}{-}$ and $\ttp{X}{-}$ preserve homotopy.
  \item A morphism $\mapdef{\a}{Y}{Y'}$ of $R$-complexes yields
    isomorphisms of $\kk$-complexes
    \begin{equation*}
      \Cone{\btpp{X}{\a}}\dis \btp{X}{\Cone{\a}}  \qand
      \Cone{\ttpp{X}{\a}}\dis \ttp{X}{\Cone{\a}}\:.\mspace{-30mu}
    \end{equation*}
    It follows that $\btp{X}{\a}$ is a quasi-isomorphism if and only
    if $\btp{X}{\Cone{\a}}$ is acyclic, and similarly for
    $\ttp{X}{\a}$.
  \end{prt}
\end{ipg}

\begin{lem}
  \label{lem:dcx}
  Let $D = (\dcx{D}{i}{j})$ be a double $\kk$-complex. Let $z$ be a
  cycle in the product totalization of $D$. Assume $z$ contains a
  component $z_{m,n}$ that satisfies $z_{m,n} = \del^h(x') +
  \del^v(x'')$ for some $x'\in D_{m+1,n}$ and $x''\in D_{m,n+1}$. If
  both
  \begin{rqm}
  \item %
    $\Hh[m+k]{\dcx{D}{\ast}{n-k}} = 0$ for every $k > 0$ \ \text{and}
  \item %
    $\Hv[n+k]{\dcx{D}{m-k}{\ast}} = 0$ for every $k > 0$,
  \end{rqm}
  then $z$ is a boundary in the product totalization of $D$.
\end{lem}

\begin{prf*}
  The goal is to prove the existence of a sequence $(x_i)$ with
  $x_{m+1}=x'$ and $x_{m}=x''$ such that
  \begin{equation*}
    z_{m+k,n-k} = \del^h(x_{m+k+1}) + \del^v(x_{m+k})
  \end{equation*}
  holds for all $k \in\ZZ$; i.e.,\ one has $z = \del(x)$ in the
  product totalization of $D$.  Assume $k$ is positive and that
  $x_{m},\ldots,x_{m+k}$ have been constructed.  Then
  \begin{equation*}
    z_{m+k-1,n-k+1} = \del^h(x_{m+k}) + \del^v(x_{m+k-1})
  \end{equation*}
  holds. As $z$ is a cycle in the product totalization of $D$, one has
  \begin{align*}
    0 &= \del^h(z_{m+k,n-k}) + \del^v(z_{m+k-1,n-k+1})\\
    &= \del^h(z_{m+k,n-k}) + \del^v(\del^h(x_{m+k}) + \del^v(x_{m+k-1}))\\
    &= \del^h(z_{m+k,n-k}) - \del^h(\del^v(x_{m+k}))\:.
  \end{align*}
  Thus, $z_{m+k,n-k} - \del^v(x_{m+k})$ is a horizontal cycle, and
  hence by \rqmlbl{1} it is a horizontal boundary. There exists,
  therefore, an element $x_{m+k+1}$ with $\del^h(x_{m+k+1}) +
  \del^v(x_{m+k}) = z_{m+k,m-k}$. By induction, this argument yields
  the elements $x_{m+k}$ for $k > 1$; a symmetric argument using
  \rqmlbl{2} yields $x_{m+k}$ for $k<0$.
\end{prf*}

\begin{prp}
  \label{prp:btp}
  Let $X$ be an $\Rop$-complex and let $Y$ be an $R$-complex.
  \begin{prt}
  \item I\fsp[.4]\ $X$ is bounded above and $\tp{X_i}{Y}$ is acyclic
    for all $i$, then $\btp{X}{Y}$ is acyclic.
  \item I\fsp\ $Y$ is bounded below and $\tp{X_i}{Y}$ is acyclic for
    all $i$, then $\btp{X}{Y}$ is acyclic.
  \end{prt}
\end{prp}

\begin{prf*}
  The product totalization of the double complex $\tpp[]{X_i}{Y_j}$ is
  $\btp{X}{Y}$.

  (a) Let $m$ be an upper bound for $X$, let $n$ be an integer, and
  let $z = (z_{m+k,n-k})_k$ be a cycle in $\btp{X}{Y}$ of degree
  $m+n$. The component $z_{m,n}$ is a cycle in $\tp{X_m}{Y}$ and hence
  a boundary as $\tp{X_m}{Y}$ is assumed to be acyclic. Moreover, the
  assumption that $\tp{X_i}{Y}$ is acyclic for every $i$ ensures that
  condition (2) in \lemref{dcx} is met. Condition (1) is satisfied due
  to the boundedness of $X$; thus $z$ is a boundary in
  $\btp{X}{Y}$. As $n$ is arbitrary, it follows that $\btp{X}{Y}$ is
  acyclic.

  (b) A similar argument with the roles of $m$ and $n$ exchanged
  handles this case.
\end{prf*}

\section{Homology}
\label{sec:homology}

\noindent
Now we consider the homology of unbounded and stable tensor product
complexes. Our notation differs slightly from the one employed by
Goichot~\cite{FGc92}.

\begin{dfn}
  \label{dfn:tTor}
  Let $M$ be an $\Rop$-module and $N$ be an $R$-module. Let $P \qra M$
  be a projective resolution and let $N \qra I$ be an injective
  resolution. The $\kk$-modules
  \begin{equation*}
    \bTor{i}{M}{N} \deq \H[i]{\btp{P}{I}}
  \end{equation*}
  are the \emph{unbounded} homology modules of $M$ and $N$ over $R$,
  and the \emph{stable} homology modules of $M$ and $N$ over $R$ are
  \begin{equation*}
    \Stor[R]{i}{M}{N} \deq \H[i+1]{\ttp{P}{I}}\:.
  \end{equation*}
\end{dfn}

\begin{ipg}
  \label{ipg:1}
  As any two projective resolutions of $M$, and similarly any two
  injective resolutions of $N$, are homotopy equivalent, it follows
  from \ref{btp}(d) that the definitions in \dfnref[]{tTor} are
  independent of the choices of resolutions.
\end{ipg}

\begin{ipg}
  \label{simplifications 2}
  Notice from \ref{simplifications} that if $M$ has finite projective
  dimension, or if $N$ has finite injective dimension, then for every
  $i \in \ZZ$ one has $\Tor{i}{M}{N} = \bTor{i}{M}{N}$ and hence
  $\Stor{i}{M}{N}=0$.
\end{ipg}

\begin{ipg}
  \label{ipg:functorial}
  The standard liftings of module homomorphisms to morphisms of
  resolutions imply that the definitions of $\bTor{i}{M}{N}$ and
  $\Stor{i}{M}{N}$ are functorial in either argument; that is, for
  every $i \in \ZZ$ there are functors
  \begin{equation*}
    \dmapdef{\bTor{i}{-}{-},\: \Stor{i}{-}{-}}
    {\Mod[\Rop]\times\Mod}{\Mod[\kk]}\:.
  \end{equation*}
  These functors are homological in the sense that every short exact
  sequence of $\Rop$-modules $0 \to M' \to M \to M'' \to 0$ and every
  $R$-module $N$ give rise to a connected exact sequence of stable
  homology modules
  \begin{equation}
    \label{eq:functorial-M}
    \cdots \to \Stor{i+1}{M''}{N} \to \Stor{i}{M'}{N} \to \Stor{i}{M}{N}
    \to \Stor{i}{M''}{N} \to \cdots,
  \end{equation}
  and to an analogous sequence of $\widebar{\operatorname{Tor}}$
  modules. Similarly, every short exact sequence $0 \to N' \to N \to
  N'' \to 0$ of $R$-modules and every $\Rop$-module $M$ give rise to a
  connected exact sequence,
  \begin{equation}
    \label{eq:functorial-N}
    \cdots \to \Stor{i+1}{M}{N''} \to \Stor{i}{M}{N'} \to \Stor{i}{M}{N}
    \to \Stor{i}{M}{N''} \to \cdots,
  \end{equation}
  and an analogous sequence of $\widebar{\operatorname{Tor}}$
  modules. For details see \seccite[1]{FGc92}.
\end{ipg}

Stable and unbounded homology entwine with absolute homology.

\begin{ipg}
  \label{ipg:2}
  For every $\Rop$-module $M$ and every $R$-module $N$ the exact
  sequence of complexes \eqref{ts} yields an exact sequence of
  homology modules:
  \begin{equation}
    \label{eq:ipg2}
    \!\cdots \to \Stor{i}{M}{N} \xra{
    } \Tor{i}{M}{N} \xra{
    }
    \bTor{i}{M}{N} \xra{
    } \Stor{i-1}{M}{N} \to \cdots\!\!\!
  \end{equation}
\end{ipg}

\subsection*{\textbf{\em Flat resolutions}} 
Just like absolute homology, unbounded and stable homology can be
computed using flat resolutions in place of projective resolutions.

\begin{prp}
  \label{prp:Stor-flat}
  Let $M$ be an $\Rop$-module and $N$ be an $R$-module. If $F \qra M$
  is a flat resolution and $N \qra I$ is an injective resolution, then
  there are isomorphisms,
  \begin{equation*}
    \bTor{i}{M}{N} \is \H[i]{\btp{F}{I}} \qand
    \Stor{i}{M}{N} \is \H[i+1]{\ttp{F}{I}}
  \end{equation*}
  for every $i\in\ZZ$, and they are natural in either argument.
\end{prp}

\begin{prf*}
  Let $P \qra M$ be a projective resolution, then there is a
  quasi-isomorphism $\mapdef{\alpha}{P}{F}$. For every $R$-complex $E$
  the complex $\tp{\Cone(\alpha)}{E}$ is acyclic by
  \corcite[7.5]{LWCHHlb}. In particular, $\tp{\Cone(\alpha)}{I}$ is
  acyclic, and it follows from \prpref{btp} that
  $\btp{\Cone(\alpha)}{I}$ is acyclic; now $\ttp{\Cone(\alpha)}{I}$ is
  acyclic in view of \eqref{ts}. Thus, the morphisms $\btp{\alpha}{I}$
  and $\ttp{\alpha}{I}$ are quasi-isomorphisms by \ref{btp}(e), and
  the desired isomorphisms follow from \dfnref[]{tTor}; it is standard
  to verify that they are natural.
\end{prf*}

In view of \ref{simplifications} the next result is now immediate.

\begin{cor}
  \label{cor:Stor-flat}
  For every $\Rop$-module $M$ of finite flat dimension,
  $\Stor{i}{M}{-}=0$ holds for every $i\in\ZZ$. \qed
\end{cor}

\begin{rmk}
  \label{rmk:pdfm}
  If every flat $\Rop$-module has finite projective dimension, then
  \corref{Stor-flat} is covered by \ref{simplifications 2}.  There are
  wide classes of rings over which flat modules have finite projective
  dimension: Iwanaga-Gorenstein rings, see \thmcite[9.1.10]{rha}, and
  more generally rings of finite finitistic projective dimension, see
  Jensen \prpcite[6]{CUJ70}; in a different direction, rings of
  cardinality at most $\aleph_n$ for some $n\in\NN$, see Gruson and
  Jensen \thmcite[7.10]{LGrCUJ81}.

  On the other hand, recall that any product of fields is a von
  Neumann regular ring, i.e., every module over such a ring is
  flat. Osofsky \cite[3.1]{BLO70} shows that a sufficiently large
  product of fields has infinite global dimension and hence must have
  flat modules of infinite projective dimension.
\end{rmk}

The next result is an analogue for stable homology of
\thmcite[2.2]{LLAOVl07}.

\begin{prp}
  \label{prp:btor0}
  Let $M$ be an $\Rop$-module and let $n\in\ZZ$. The following
  conditions are equivalent.
  \begin{eqc}
  \item The connecting morphism $\Stor{i}{M}{-} \to \Tor{i}{M}{-}$ is
    an isomorphism for every $i\ge n$.
  \item $\Tor{i}{M}{E}=0$ for every injective $R$-module $E$ and every
    $i\ge n$.
  \item $\bTor{i}{M}{-}=0$ for every $i\ge n$.
  \end{eqc}
\end{prp}

\begin{prf*}
  For every injective $R$-module $E$ and for every $i\in\ZZ$ one has
  $\Stor{i}{M}{E}=0$ by \ref{simplifications 2}. Thus \eqclbl{i}
  implies \eqclbl{ii}.

  In the following we use the notation $\Tsb{s}{(-)}$ for the soft
  truncation below at $s$.

  \proofofimp{ii}{iii} Let $N$ be an $R$-module and $N \qra I$ be an
  injective resolution.  Let $P \qra M$ be a projective resolution and
  let $P^+$ denote its mapping cone. Let $E$ be an injective
  $R$-module; by assumption one has $\Tor{i}{M}{E}=0$ for all $i\ge
  n$. By right exactness of the tensor product,
  $\tp{\Tsb{n-2}{P^+}}{E}$ is acyclic and hence
  $\btp{\Tsb{n-2}{P^+}}{I}$ is acyclic by \ref{btp}(a) and
  \prpref{btp}. Thus, for every $i\ge n$ one has
  \begin{equation*}
    \bTor{i}{M}{N} = \H[i]{\btp{P}{I}}
    = \H[i]{\btp{\Tsb{n-2}{P^+}}{I}} = 0.
  \end{equation*}
  (Notice that for all $n\le 0$ one has $\Tsb{n-2}{P^+} = P^+$.)

  \proofofimp{iii}{i} This is immediate from the exact sequence
  \eqref{ipg2}.
\end{prf*}

As vanishing of $\Tor{\gg 0}{M}{-}$ detects finite flat dimension of
$M$, the next result is an immediate consequence of \prpref{btor0}. In
the parlance of \cite{EEnOJn93a} it says that the flat and copure flat
dimensions agree for $\Rop$-modules $M$ with $\Stor{}{M}{-}=0$.

\begin{cor}
  \protect\pushQED{\qed}%
  \label{cor:copureflat}
  Let $M$ be an $\Rop$-module. If\, $\Stor{i}{M}{-}=0$ for all
  $i\in\ZZ$, then
  \begin{equation*}
    \fd[\Rop]{M} = \sup\setof{i\in\ZZ}{\Tor{i}{M}{E}\ne0 \text{ for
        some injective
        $R$-module }E}\:.\qedhere
  \end{equation*}
\end{cor}

Dimension shifting is a useful tool in computations with stable
homology.

\begin{ipg}
  \label{ipg:(co)syz}
  For a module $X$, we denote by $\Omega_{m}X$ the $m^{\rm th}$ syzygy
  in a projective resolution of $X$ and by $\Omega^{m}X$ the $m^{\rm
    th}$ cosyzygy in an injective resolution of $X$. By Schanuel's
  lemma a syzygy/cosyzygy is uniquely determined up to a
  projective/injective summand. In view of \ref{simplifications 2} and
  the exact sequences in \ref{ipg:functorial} one gets
  \begin{alignat}{2}
    \label{eq:dimshiftproj}
    \Stor{i}{M}{N} &\dis \Stor{i-m}{\Omega_mM}{N}&&\quad\text{for } m\ge 0\qand\\
    \label{eq:dimshiftinj}
    \Stor{i}{M}{N} &\dis \Stor{i+m}{M}{\Omega^{m}N}&& \quad\text{for }
    m\ge 0\:.
  \end{alignat}
  Moreover, if $F \qra M$ is a flat resolution of $M$, then there are
  isomorphisms,
  \begin{equation}
    \label{eq:dimshiftflat}
    \Stor{i}{M}{N} \dis \Stor{i-m}{\Co[m]{F}}{N} \quad\text{for } m\ge 0\:;\\
  \end{equation}
  this follows from \eqref{functorial-M} and \corref{Stor-flat}.
\end{ipg}

\begin{ipg}
  \label{S-structure}
  Suppose $R$ is commutative and $S$ is a flat $R$-algebra.  For every
  $\Sop$-module $M$, a projective resolution $P \qra M$ over $\Sop$ is
  a flat resolution of $M$ as an $R$-module. Thus, it follows from
  \prpref{Stor-flat} that stable homology $\Stor{i}{M}{-}$ is a
  functor from $\Mod[R]$ to $\Mod[\Sop]$. Similarly, an injective
  resolution $N \qra I$ of an $S$-module is also an injective
  resolution of $N$ as an $R$-module, whence $\Stor{i}{-}{N}$ is a
  functor from $\Mod$ to $\Mod[S]$.
\end{ipg}

\section{Vanishing of stable homology
  $\widetilde{\mathrm{Tor}}\mspace{1mu}(M,-)$}
\label{sec:vanishing1}

\noindent
We open with a partial converse to \corref{Stor-flat}; recall, for
example from \prpcite[3.2.12]{rha}, that a finitely presented module
is flat if and only if it is projective.

\begin{thm}
  \label{thm:finitepd}
  Let $M$ be an $\Rop$-module that has a degree-wise finitely
  generated projective resolution. The following conditions are
  equivalent.
  \begin{eqc}
  \item $\pd[\Rop]{M}$ is finite.
  \item $\Stor{i}{M}{-}=0$ for all $i\in\ZZ$.
  \item $\Stor{i}{M}{-}=0$ for some $i \ge 0$.
  \item $\Stor{0}{M}{\Hom[\kk]{M}{E}}=0$ for some faithfully injective
    $\kk$-module $E$.
  \end{eqc}
  Moreover, if $R$ is left noetherian with $\id{R}$ finite, then
  \eqclbl{i}--\eqclbl{iv} are equivalent to
  \begin{eqc}
  \item[\eqclbl{iii'}] $\Stor{i}{M}{-}=0$ for some $i\in\ZZ$.
  \end{eqc}
\end{thm}

\begin{prf*}
  The implications
  \eqclbl{i}$\implies$\eqclbl{ii}$\implies$\eqclbl{iii}$\implies$\eqclbl{iii'}
  are clear; see \ref{simplifications 2}. Part \eqclbl{iv} follows
  from \eqclbl{iii} by dimension shifting \eqref{dimshiftinj}.

  For \proofofimp[]{iv}{i} use \thmref{stable-Ext-Tor-1} to get
  \begin{equation*}
    0 = \Stor{0}{M}{\Hom[\kk]{M}{E}} \dis \Hom[\kk]{\Sext[\Rop]{0}{M}{M}}{E}\:.
  \end{equation*}
  Thus $\Sext[\Rop]{0}{M}{M}$ is zero, hence $M$ has finite projective
  dimension by \prpcite[2.2]{LLAOVl07}.

  Finally, assume that $R$ is left noetherian with $\id{R}$
  finite. Let $E$ be a faithfully injective $\kk$-module; in view of
  what has already been proved, it is sufficient to show that
  $\Stor{-i}{M}{-}=0$ for some $i > 0$ implies
  $\Stor{0}{M}{\Hom[\kk]{M}{E}}=0$. It follows from the assumptions on
  $R$ that every free $R$-module has finite injective dimension. Thus
  every projective $R$-module $P$ has finite injective dimension, and
  \pgref{simplifications 2} yields $\Stor{i}{-}{P}=0$ for all
  $i\in\ZZ$. From the exact sequence \eqref{functorial-N} it follows
  that there are isomorphisms $\Stor{0}{M}{\Hom[\kk]{M}{E}} \is
  \Stor{-i}{M}{\Omega_{i}\Hom[\kk]{M}{E}}$ for $i > 0$, so vanishing
  of $\Stor{-i}{M}{-}$ implies $\Stor{0}{M}{\Hom[\kk]{M}{E}}=0$.
\end{prf*}

\begin{cor}
  \label{cor:artin pd}
  Let $R$ be an Artin algebra with duality functor $\Dual{-}$. For a
  finitely generated $\Rop$-module $M$ the following conditions are
  equivalent.
  \begin{eqc}
  \item $\pd[\Rop]{M}$ is finite.
  \item There is an $i \ge 0$ with $\Stor{i}{M}{N}=0$ for all finitely
    generated $R$-modules~$N$.
  \item $\Stor{0}{M}{\Dual{M}}=0$.
  \end{eqc}
\end{cor}

\begin{prf*}
  Part \eqclbl{ii} follows from \eqclbl{i} in view of
  \ref{simplifications 2}. Part \eqclbl{iii} follows from \eqclbl{ii}
  by dimension shifting \eqref{dimshiftinj} in a degree-wise finitely
  generated injective resolution of $\Dual{M}$. Finally, one has
  $\Dual{M} = \Hom[\kk]{M}{E}$, see the first paragraph in
  \secref{tensor}, so \eqclbl{i} follows
  from \eqclbl{iii} by \thmref{finitepd}, as the injective hull $E$ of
  $\kk/\mathrm{Jac}(\kk)$ is faithfully injective.
\end{prf*}

By a result of Bass and Murthy \lemcite[4.5]{HBsMPM67}, a commutative
noetherian ring is regular if and only if every finitely generated
module has finite projective dimension.

\begin{cor}
  \label{cor:regular}
  For a commutative noetherian ring $R$, the following are equivalent:
  \begin{eqc}
  \item $R$ is regular.
  \item $\Stor{i}{M}{-}=0$ for all finitely generated $R$-modules $M$
    and all $i\in\ZZ$.
  \item There is an $i\in\ZZ$ with $\Stor{i}{M}{-}=0$ for all finitely
    generated $R$-modules~$M$.
  \end{eqc}
\end{cor}

\begin{prf*}
  The implications
  \eqclbl{i}$\implies$\eqclbl{ii}$\implies$\eqclbl{iii} are clear; see
  \ref{simplifications 2}.

  \proofofimp{iii}{i} If $i\ge 0$, then every finitely generated
  $R$-module has finite projective dimension by \thmref{finitepd};
  i.e., $R$ is regular. If $i<0$, then for any finitely generated
  $R$-module $M$ one has $\Stor{0}{M}{-} \is \Stor{i}{\Omega_{-i}M}{-}
  =0$ by dimension shifting \eqref{dimshiftproj}, and as above it
  follows that $R$ is regular.
\end{prf*}

\subsection*{\textbf{\em Tate flat resolutions}} 
We show that under extra assumptions on $M$ one can compute
$\Stor{}{M}{N}$ without resolving the second argument.

\begin{ipg}
  \label{ipg:C}
  Let $N$ be an $R$-module; choose an injective resolution \mbox{$N
    \qra I$} and a projective resolution $P \qra N$. The composite $P
  \to N \to I$ is a quasi-isomorphism, so the complex $\cone{N}{PI} =
  \Cone{(P \to N \to I)}$ is acyclic, and up to homotopy equivalence
  it is independent of the choice of resolutions. This construction is
  functorial in $N$.
\end{ipg}

\begin{lem}
  \label{lem:tTor}
  Let $F$ be a bounded below complex of flat $\Rop$-modules; let $N$
  be an $R$-module with a projective resolution $P \qra N$ and an
  injective resolution \mbox{$N\qra I$}.  There is an isomorphism
  $\ttp{F}{I} \qis \btp{F}{\cone{N}{PI}}$ in the derived
  category~$\D[\kk]$, and it is functorial in the second argument;
  cf.~\ref{ipg:C}.
\end{lem}

\begin{prf*}
  There is by \ref{btp}(c) a short exact sequences of $\kk$-complexes
  \begin{equation*}
    0 \lra \btp{F}{I} \lra
    \btp{F}{\cone{N}{PI}} \lra
    \btp{F}{\susp P} \lra 0\:.
  \end{equation*}
  With \eqref{ts} it forms a commutative diagram whose rows are
  triangles in~$\D[\kk]$,
  \begin{equation*}
    \xymatrix{
      \btp{F}{P} \ar[r]
      \ar[d]_-{\qis}
      & \btp{F}{I} \ar[r] \ar@{=}[d]
      & \btp{F}{\cone{N}{PI}} \ar[r]
      & \susp\btpp{F}{P} \ar[d]_-{\qis}\\
      \tp{F}{I} \ar[r]
      & \btp{F}{I} \ar[r]
      & \ttp{F}{I} \ar[r]
      & \susp \tpp{F}{I}\:,
    }
  \end{equation*}
  where the top triangle is rotated back once. The isomorphism is
  induced by the composite quasi-isomorphism $P \to N \to I$, as one
  has $\btp{F}{P} = \tp{F}{P}$; cf.~\ref{simplifications}. Completion
  of the morphism of triangles yields the desired isomorphism via the
  Triangulated Five Lemma; see \prpcite[4.3]{THlPJr10}. Functoriality
  in the second argument is straightforward to verify.
\end{prf*}

\begin{lem}
  \label{lem:N}
  Let $N$ be an $R$-module with a projective resolution $P \qra N$ and
  an injective resolution $N \qra I$. Let $T$ be a complex of flat
  $\Rop$-modules such that $\tp{T}{E}$ is acyclic for every injective
  $R$-module $E$. There is then an isomorphism $\btp{T}{\cone{N}{PI}}
  \qis \susp\tpp{T}{N}$ in $\D[\kk]$ which is functorial in the second
  argument; cf.~\ref{ipg:C}.
\end{lem}

\begin{prf*}
  By the assumptions and \partprpref{btp}{a} the complex $\btp{T}{I}$
  is acyclic, so application of $\btp{T}{-}$ to the exact sequence $0
  \to I \to \cone{N}{PI} \to \susp P \to 0$ yields a quasi-isomorphism
  $\btp{T}{\cone{N}{PI}} \to \btp{T}{\susp P}$; cf.~\ref{btp}(c).
  Now, let $P^+$ denote the mapping cone of the quasi-isomorphism $P
  \to N$. The complex $\btp{T}{P^+}$ is acyclic
  by \partprpref{btp}{b}, so application of $\btp{T}{-}$ to the
  mapping cone sequence $0 \to N \to P^+ \to \susp P \to 0$ yields an
  isomorphism $\btp{T}{\susp P} \qis \susp(\btp{T}{N})$ in $\D[\kk]$,
  and the desired isomorphism follows as $N$ is a module;
  cf.~\ref{simplifications}. Functoriality in the second argument is
  straightforward to verify; cf.~\ref{ipg:C}.
\end{prf*}

\begin{rmk}
  The quasi-isomorphisms in \lemref[Lemmas~]{tTor} and \lemref[]{N}
  hold with $\cone{N}{PI}$ replaced by $\cone{N}{FI}$ where $F \qra N$
  is a flat resolution.
\end{rmk}

The objects discussed in the next paragraph appear in the literature
under a variety of names; here we stick with the terminology from
\cite{LLn13}.

\begin{ipg}
  \label{ipg:tateflat}
  An acyclic complex $T$ of flat $\Rop$-modules is called
  \emph{totally acyclic} if $\tp{T}{E}$ is acyclic for every injective
  $R$-module $E$.  A \emph{Tate flat resolution} of an $\Rop$-module
  $M$ is a pair $(T,F)$, where $F\qra M$ is a flat resolution and $T$
  is a totally acyclic complex of flat $\Rop$-modules with $\Thb{n}{T}
  \is \Thb{n}{F}$ for some $n\ge 0$. Tate's name is invoked here
  because these resolutions can be used to compute Tate homology; see
  \thmcite[A]{LLn13}.
\end{ipg}

\begin{ipg}
  \label{ipg:gflatTate}
  Over a left coherent ring $R$, an $\Rop$-module $M$ has a Tate flat
  resolution if and only if it has finite Gorenstein flat dimension;
  see \prpcite[3.4]{LLn13}.

  If $R$ is noetherian, then every finitely generated $\Rop$-module
  $M$ of finite G-dimension has a Tate flat resolution. Here is a
  direct argument: By \thmcite[3.1]{LLAAMr02} there exists a pair
  $(T,P)$ where $P \qra M$ is a projective resolution and $T$ is an
  acyclic complex of finitely generated projective $\Rop$-modules with
  the following properties: the complex $\Hom[\Rop]{T}{R}$ is acyclic,
  and there is an $n\ge 0$ with $\Thb{n}{T} \is \Thb{n}{P}$. For every
  injective $R$-module $E$ one now has
  \begin{equation}
    \label{eq:gflatTate}
    \tp{T}{E} \is \tp{T}{\Hom{R}{E}} \is \Hom{\Hom[\Rop]{T}{R}}{E}
  \end{equation}
  where the last isomorphism holds because $T$ is degree-wise finitely
  generated; see~\prpcite[VI.5.2]{careil}. Thus $T$ is a totally
  acyclic complex of flat $\Rop$-modules and $(T,P)$ is a Tate flat
  resolution of $M$.
\end{ipg}

The next result shows that stable homology can be computed via Tate
flat resolutions. We apply it in \secref{tate}
to compare stable homology to Tate homology, but before that it is
used in the proofs of \thmref[Theorems~]{vanishing} and
\thmref[]{Stor-balanced}.

\begin{thm}
  \label{thm:Stor-Tateflat}
  Let $M$ be an $\Rop$-module that has a Tate flat resolution
  $(T,F)$. For every $R$-module $N$ and every $i\in\ZZ$ there is an
  isomorphism,
  \begin{equation*}
    \Stor{i}{M}{N} \dis \H[i]{\tp{T}{N}}\:,
  \end{equation*}
  and it is functorial in the second argument.
\end{thm}

\begin{prf*}
  Choose $n\ge 0$ such that $\Thb{n}{T}$ and $\Thb{n}{F}$ are
  isomorphic, and consider the degree-wise split exact sequence
  \begin{equation*}
    0 \lra \Tha{n-1}{T} \lra T \lra \Thb{n}{F} \lra 0\:.
  \end{equation*}
  Let $P \qra N$ be a projective resolution and $N \qra I$ be an
  injective resolution; set $C = \cone{N}{PI}$; cf.~\ref{ipg:C}. There
  is, by \ref{btp}(c), an exact sequence of $\kk$-complexes,
  \begin{equation*}
    0 \lra \btp{\Tha{n-1}{T}}{C} \lra \btp{T}{C} \lra \btp{\Thb{n}{F}}{C}
    \lra 0\:.
  \end{equation*}
  As $C$ is acyclic, it follows from \partprpref{btp}{a} that
  $\btp{\Tha{n-1}{T}}{C}$ is acyclic, whence the surjective morphism
  above is a quasi-isomorphism. Since $T$ is a totally acyclic complex
  of flat modules, there is by \lemref{N} an isomorphism $\btp{T}{C}
  \qis \susp\tpp{T}{N}$ in $\D[\kk]$.  Moreover, \lemref{tTor} yields
  an isomorphism $\btp{\Thb{n}{F}}{C} \qis \ttp{\Thb{n}{F}}{I}$ in
  $\D[\kk]$. Thus one has $\susp \tpp{T}{N} \qis \ttp{\Thb{n}{F}}{I}$
  in $\D[\kk]$. This explains the second isomorphism in the next
  computation.
  \begin{align*}
    \Stor{i}{M}{N}
    &\dis \Stor{i-n}{\Co[n]{F}}{N}\\
    &\deq \H[i-n+1]{\ttp{(\susp^{-n}\Thb{n}{F})}{I}}\\
    &\deq \H[i+1]{\ttp{\Thb{n}{F}}{I}}\\
    &\dis \H[i+1]{\susp \tpp{T}{N}}\\
    &\deq \H[i]{\tp{T}{N}}
  \end{align*}
  The first isomorphism holds by \eqref{dimshiftflat}; the first
  equality follows from \prpref{Stor-flat} as the canonical surjection
  $\susp^{-n}\Thb{n}{F} \to \Co[n]{F}$ is a flat resolution.
\end{prf*}

As discussed in the introduction, the G-dimension is a homological
invariant for finitely generated modules over noetherian
rings. Characterizations of modules of finite G-dimension have
traditionally involved both vanishing of (co)homology and
invertibility of a certain morphism; see for example \cite[(2.1.6),
(2.2.3), (3.1.5), and (3.1.11)]{lnm}. More recently, Avramov, Iyengar,
and Lipman \cite{AIL-10} showed that a finitely generated module $M$
over a commutative noetherian ring $R$ has finite G-dimension if and
only $M$ is isomorphic to the complex $\RHom{\RHom{M}{R}}{R}$ in the
derive category $\D$. The crucial step in our proof of the next
theorem is to show that vanishing of stable homology $\Stor{}{M}{R}$
implies such an isomorphism.

\begin{thm}
  \label{thm:vanishing}
  Let $R$ be a commutative noetherian ring. For a finitely generated
  $R$-module $M$, the following conditions are equivalent.
  \begin{eqc}
  \item $\Gdim{M}<\infty$.
  \item $\Stor[R]{i}{M}{N}=0$ for every $R$-module $N$ of finite flat
    dimension and all $i\in \ZZ$.
  \item $\Stor[R]{i}{M}{R}=0$ for all $i\in \ZZ$.
  \end{eqc}
\end{thm}

\begin{prf*}
  The implication $\proofofimp[]{ii}{iii}$ is trivial.

  \proofofimp{i}{ii} By \pgref{ipg:gflatTate} the $\Rop$ module $M$
  has a Tate flat resolution $(T,F)$. \thmref{Stor-Tateflat} yields
  isomorphisms $\Stor[R]{i}{M}{-} \is \H[i]{\tp{T}{-}}$ for all
  $i\in\ZZ$. It follows by induction on the flat dimension of $N$ that
  $\tp{T}{N}$ is acyclic; cf.~\lemcite[2.3]{CFH-06}.

  \proofofimp{iii}{i} Let $P \qra M$ be a degree-wise finitely
  generated projective resolution, and let $R \qra I$ be an injective
  resolution. By assumption, the complex $\ttp{P}{I}$ is acyclic,
  which explains the third $\qis$ below:
  \begin{equation*}
    M \dqis \tp{P}{R}
    \dqis\tp{P}{I}
    \dqis\btp{P}{I}
    \dis\btp{P}{\Hom{R}{I}}
    \dis\Hom{\Hom{P}{R}}{I}\:.
  \end{equation*}
  The last isomorphism holds by \prpref{bhev} as $\bHom{P}{R}$ equals
  $\Hom{P}{R}$. Now it follows from \thmcite[2]{AIL-10} and
  \thmcite[(2.2.3)]{lnm} that $\Gdim{M}$ is finite.
\end{prf*}

\begin{cor}
  \label{cor:Gor1}
  For a commutative noetherian ring $R$ the following are equivalent.
  \begin{eqc}
  \item $R$ is Gorenstein.
  \item $\Stor{i}{M}{R}=0$ for all finitely generated $R$-modules $M$
    and all $i\in\ZZ$.
  \item $\Stor{i}{M}{R}=0$ for all finitely generated $R$-modules $M$
    and all $i< 0$.
  \end{eqc}
\end{cor}

\begin{prf*}
  By a result of Goto, $R$ is Gorenstein if and only if every finitely
  generated $R$-module has finite G-dimension
  \corcite[2]{SGt82}. Combining this with \thmref{vanishing} one gets
  the the equivalence of \eqclbl{i} and \eqclbl{ii}. The implication
  \eqclbl{ii}$\implies$\eqclbl{iii} is clear. Finally, \eqclbl{ii}
  follows from \eqclbl{iii} by dimension shifting
  \eqref{dimshiftproj}.
\end{prf*}

\section{Balancedness of stable homology}
\label{sec:balance}

\noindent
Absolute homology is balanced over any ring: there are always
isomorphisms $\Tor{}{M}{N} \is \Tor[\Rop]{}{N}{M}$. It follows already
from \corref{artin pd} that stable homology can be balanced only over
special rings. Indeed, if $R$ is an Artin algebra and stable homology
is balanced over $R$, then the dual module of $R$ has finite
projective dimension over both $R$ and $\Rop$, whence $R$ is
Iwanaga-Gorenstein. The converse is part of \corref{art-balance}.

We open this section with a technical lemma; it is similar to a result
of Enochs, Estrada, and Iacob \thmcite[3.6]{EEI-12}. Recall the
notation $\Co[m]{-}$ from \ref{intro notation}.

\begin{lem}
  \label{lem:complex-balanced}
  Let $T$ and $T'$ be acyclic complexes of flat $\Rop$-modules and
  flat $R$-modules, respectively. For all integers $m$ and $n$ there
  are isomorphisms in $\D[\kk]$,
  \begin{align*}
    \tp{\Tha{n-1}{T}}{\Co[m]{T'}} &\qis \susp^{n-m}
    (\tp{\Co[n]{T}}{\Tha{m-1}{T'}})\\
    \tp{\Thb{n}{T}}{\Co[m]{T'}} &\qis \susp^{n-m}
    (\tp{\Co[n]{T}}{\Thb{m}{T'}})\\
    \tp{T}{\Co[m]{T'}} &\qis \susp^{n-m} (\tp{\Co[n]{T}}{T'})\:.
  \end{align*}
\end{lem}

\begin{prf*}
  Because the complexes $T$ and $T'$ are acyclic, there are
  quasi-isomorphisms $\Co[n]{T} \to \susp^{1-n}\Tha{n-1}{T}$ and
  $\Co[m]{T'} \to \susp^{1-m}\Tha{m-1}{T'}$. Hence
  \prpcite[2.14(b)]{CFH-06} yields
  \begin{align*}
    \tp{\Tha{n-1}{T}}{\Co[m]{T'}}
    &\qis \tp{\Tha{n-1}{T}}{\susp^{1-m}\Tha{m-1}{T'}}\\
    &\qis \susp^{n-m}(\tp{\susp^{1-n}\Tha{n-1}{T}}{\Tha{m-1}{T'}})\\
    &\qis \susp^{n-m} (\tp{\Co[n]{T}}{\Tha{m-1}{T'}})\:.
  \end{align*}
  This demonstrates the first isomorphism in the statement, and the
  second is proved similarly.  Finally, these two isomorphisms connect
  the exact sequences
  \begin{equation*}
    0 \to \tp{\Tha{n-1}{T}}{\Co[m]{T'}} \to \tp{T}{\Co[m]{T'}} \to
    \tp{\Thb{n}{T}}{\Co[m]{T'}} \to 0
  \end{equation*}
  and
  \begin{equation*}
    0 \to \tp{\Co[n]{T}} {\Tha{m-1}{T'}}\to \tp{\Co[n]{T}}{T'}
    \to \tp{\Co[n]{T}}{\Thb{m}{T'}} \to 0\:,
  \end{equation*}
  and one obtains the last isomorphism via the Triangulated Five Lemma in
  $\D[\kk]$; see \prpcite[4.3]{THlPJr10}.
\end{prf*}

\begin{thm}
  \label{thm:Stor-balanced}
  Let $M$ be an $\Rop$-module and $N$ be an $R$-module. If they both
  have Tate flat resolutions, then for each $i\in \ZZ$ there is an
  isomorphism
  \begin{equation*}
    \Stor{i}{M}{N} \dis \Stor[\Rop]{i}{N}{M}\:.
  \end{equation*}
\end{thm}

\begin{prf*}
  Let $(T,F)$ and $(T',F')$ be Tate flat resolutions of $M$ and $N$,
  respectively. Choose $n\in\ZZ$ such that there are isomorphisms
  $\Thb{n}{T} \is \Thb{n}{F}$ and $\Thb{n}{T'} \is \Thb{n}{F'}$. For
  every $i\in\ZZ$ one has
  \begin{align*}
    \Stor{i}{M}{N}
    &\is \H[i]{T\otimes_{R}N}\\
    &\is \H[i-n]{\tp{T}{\Co[n]{T'}}}\\
    &\is \H[i-n]{\tp{\Co[n]{T}}{T'}}\\
    &\is \H[i]{\tp{M}{T'}}\\
    &\is \Stor[\Rop]{i}{N}{M}\:,
  \end{align*}
  where the first and the last isomorphisms follow from
  \thmref{Stor-Tateflat}, the second and fourth isomorphisms follow by
  dimension shifting, and the third isomorphism holds by the last
  isomorphism in \lemref{complex-balanced}.
\end{prf*}

\begin{dfn}
  Let $M$ be an $\Rop$-module and $N$ be an $R$-module.  Stable
  homology is \emph{balanced} for $M$ and $N$ if one has
  $\Stor{i}{M}{N} \is \Stor[\Rop]{i}{N}{M}$ for all $i\in\ZZ$.
\end{dfn}

\begin{ipg}
  \label{ipg:IGor}
  \thmref{Stor-balanced} says that stable homology is balanced for all
  (pairs of) $\Rop$- and $R$-modules that have Tate flat
  resolutions. If $R$ is Iwanaga--Gorenstein, then every $\Rop$-module
  and every $R$-module has a Tate flat resolution, see
  \pgref{ipg:gflatTate} and \thmcite[12.3.1]{rha}, so stable homology
  is balanced for all (pairs of) $\Rop$- and $R$-modules.
\end{ipg}

\begin{cor}
  \label{cor:art-balance}
  For an Artin algebra $R$ the following conditions are equivalent.
  \begin{eqc}
  \item $R$ is Iwanaga-Gorenstein.
  \item Stable homology is balanced for all $\Rop$- and $R$-modules.
  \item Stable homology is balanced for all finitely generated $\Rop$-
    and $R$-modules.
  \end{eqc}
\end{cor}

\begin{prf*}
  Per \pgref{ipg:IGor}, part \eqclbl{i} implies \eqclbl{ii}, which
  clearly implies \eqclbl{iii}. Let $E=\Dual{R}$ be the dual module of
  $R$; it is injective and finitely generated over $\Rop$ and over
  $R$. Thus, if stable homology is balanced for finitely generated
  $\Rop$- and $R$ modules, then it follows from \ref{simplifications
    2} and \corref{artin pd} that $\pd{E}$ as well as $\pd[\Rop]{E}$
  is finite. By duality, both $\id{R}$ and $\id[\Rop]{R}$ are then
  finite.
\end{prf*}

\begin{cor}
  \label{cor:Gor2}
  A commutative noetherian ring $R$ is Gorenstein if and only if
  stable homology is balanced for all finitely generated $R$-modules.
\end{cor}

\begin{prf*}
  Over a Gorenstein ring, the G-dimension of every finitely generated
  module is finite by \corcite[2]{SGt82}, so every finitely generated
  module has a Tate flat resolution; see \pgref{ipg:gflatTate}. Thus,
  balancedness of stable homology follows from
  \thmref{Stor-balanced}. Conversely, balancedness of stable homology
  implies $\Stor{i}{M}{R}=0$ for every finitely generated $R$-module
  $M$ and all $i\in\ZZ$, so $R$ is Gorenstein by \corref{Gor1}.
\end{prf*}

\begin{cor}
  \label{cor:Gor3}
  A commutative noetherian ring $R$ of finite Krull dimension is
  Gorenstein if and only if stable homology is balanced over $R$.
\end{cor}

\begin{prf*}
  If $R$ is Gorenstein of finite Krull dimension, then it is
  Iwanaga-Gorenstein. Therefore, stable homology is balanced over $R$ per
  \pgref{ipg:IGor}. The converse holds by \corref{Gor2}.
\end{prf*}

\section{Vanishing of stable homology
  $\widetilde{\mathrm{Tor}}\mspace{1mu}(-,N)$}
\label{sec:vanishing2}

\noindent
Vanishing of stable homology $\Stor{i}{-}{N}$ over an Artin algebra
can by duality be understood via vanishing of $\Stor{i}{M}{-}$.

\begin{prp}
  \label{prp:artin id}
  Let $R$ be an Artin algebra with duality functor $\Dual{-}$. For a
  finitely generated $R$-module $N$ the following conditions are
  equivalent.
  \begin{eqc}
  \item $\id{N}$ is finite.
  \item $\Stor{i}{-}{N}=0$ for all $i\in\ZZ$.
  \item There is an integer $i \le 0$ with $\Stor{i}{M}{N}=0$ for all
    finitely generated $\Rop$-modules~$M$.
  \item $\Stor{0}{\Dual{N}}{N}=0$.
  \end{eqc}
\end{prp}

\begin{prf*}
  The implications
  \eqclbl{i}$\implies$\eqclbl{ii}$\implies$\eqclbl{iii} are clear; see
  \ref{simplifications 2}.  Part \eqclbl{iv} follows from \eqclbl{iii}
  by dimension shifting \eqref{dimshiftproj}, as $\Dual{N}$ is
  finitely generated. Finally, vanishing of $\Stor{0}{\Dual{N}}{N} \is
  \Stor{0}{\Dual{N}}{\Dual{\Dual{N}}}$ implies by \corref{artin pd}
  that $\pd[\Rop]{\Dual{N}}$ is finite, whence $\id{N}$ is finite, and
  so \eqclbl{i} follows from \eqclbl{iv}.
\end{prf*}

\subsection*{\textbf{\em Local rings}} 
To analyze vanishing of stable homology $\Stor{i}{-}{N}$ over
commutative noetherian rings, we start locally.

\begin{ipg}
  \label{ipg:locally}
  Let $R$ be a commutative noetherian local ring with residue
  field~$k$. For an $R$-module $M$, the \emph{depth} invariant can be
  defined as
  \begin{equation*}
    \dpt{M} = \inf\setof{i\in\ZZ}{\Ext{i}{k}{M} \ne 0}\:,
  \end{equation*}
  and if $M$ is finitely generated, then its injective dimension can
  be computed as
  \begin{equation*}
    \id{M} = \sup\setof{i\in\ZZ}{\Ext{i}{k}{M} \ne 0}\:.
  \end{equation*}
  The depth is finite for $M\ne 0$. The ring $R$ is Cohen--Macaulay if
  there exists a finitely generated module $M\ne 0$ of finite
  injective dimension; this is a consequence of the New Intersection
  Theorem due to Peskine and Szpiro \cite{CPsLSz73} and
  Roberts~\cite{PRb87}.
\end{ipg}

The following is an analogue of \thmcite[6.1]{LLAOVl07}.

\begin{lem}
  \label{lem:finiteness inj dims}
  Let $R$ be a commutative noetherian local ring with residue field
  $k$. For every finitely generated $R$-module $N$ and for every
  $i\in\ZZ$ there is an isomorphism
  \begin{equation*}
    \bTor{i}{k}{N} \dis \prod_{j\in\ZZ}
    \Hom[k]{\Ext{j}{k}{R}}{\Ext{j-i}{k}{N}}
  \end{equation*}
  of $R$-modules, and in particular, $\bTor{i}{k}{N}$ is a $k$-vector
  space.
\end{lem}

\begin{prf*}
  Let $P \qra k$ be a degree-wise finitely generated projective
  resolution and let $N \qra I$ and $R \qra J$ be injective
  resolutions.  By definition $\bTor{i}{k}{N}$ is the $i^{\rm th}$
  homology of the complex %
  $\btp{P}{I}$ which we compute using \prpref{bhev} as follows
  $\btp{P}{I} \dis \btp{P}{\Hom{R}{I}} \dis \Hom{\Hom{P}{R}}{I}$. Next
  we simplify using quasi-isomorphisms and Hom-tensor adjointness:
  \begin{align*}
    \Hom{\Hom{P}{R}}{I}  & \dqis \Hom{\Hom{P}{J}}{I} \\
    & \dqis \Hom{\Hom{k}{J}}{I} \\
    & \dis \Hom{\tp[k]{\Hom{k}{J}}{k}}{I} \\
    & \dis \Hom[k]{\Hom{k}{J}}{\Hom{k}{I}}\:.
  \end{align*}
  Finally, pass to homology.
\end{prf*}


\begin{rmk}
  \lemref{finiteness inj dims} suggests that stable homology
  $\Stor{i}{k}{N}$ may be a $k$-vector space, and
  indeed it is. If $R$ is a ring and if $x$ annihilates the
  $\Rop$-module $M$, or the $R$-module $N$, then there are
  two homotopic lifts---zero and multiplication by $x$---to the
  projective resolution in the case of $M$ or to the injective
  resolution in the case of $N$; see Definition 2.1 and \cite[22.6
  and 2.3.7]{Wei}. Hence multiplication by $x$ is zero on unbounded
  and on stable homology of $M$ against $N$.  In particular, stable
  homology $\Stor{i}{k}{N}$ is a $k$-vector space.
\end{rmk}

\begin{prp}
  \label{prp:fg}
  Let $R$ be a commutative noetherian local ring with residue field
  $k$ and let $N$ be a finitely generated $R$-module. If for some
  $i\in\ZZ$ the $k$-vector space $\Stor{i}{k}{N}$ has finite rank,
  then $N$ has finite injective dimension, or $R$~is~Gorenstein.
\end{prp}

\begin{prf*}
  Each $k$-vector space $\Tor{j}{k}{N}$ has finite rank, so
  $\bTor{i+1}{k}{N}$ has finite rank by the assumption and the exact
  sequence \eqref{ipg2}.  For a finitely generated $R$-module $M\ne 0$
  the vector spaces $\Ext{j}{k}{M}$ are non-zero for all $j$ between
  $\dpt{M}< \infty$ and $\id{M}$; see Roberts~\thmcite[2]{PRb76}. When
  $\bTor{i+1}{k}{N}$ has finite rank, it follows from
  \lemref{finiteness inj dims} that $R$ or $N$ has finite injective
  dimension.
\end{prf*}

Compared to the characterization of globally Gorenstein rings in
\corref{Gor1}, condition \eqclbl{iii} below is sharper.

\begin{thm}
  \label{thm:gorenstein}
  Let $R$ be a commutative noetherian local ring with residue field
  $k$. The following conditions are equivalent.
  \begin{eqc}
  \item $R$ is Gorenstein.
  \item $\Stor[R]{i}{-}{R}=0$ for all $i\in\ZZ$.
  \item $\Stor[R]{i}{k}{R}$ has finite rank for some $i\in\ZZ$.
  \item There exists a finitely generated $R$-module $M$ such that
    $\Stor[R]{i}{k}{M}$ has finite rank for some $i\in\ZZ$ and
    $\Stor[R]{i}{M}{R}=0$ holds for all $i\in\ZZ$.
  \end{eqc}
\end{thm}

\begin{prf*}
  The implications
  \eqclbl{i}$\implies$\eqclbl{ii}$\implies$\eqclbl{iii}$\implies$\eqclbl{iv}
  are clear; see \ref{simplifications 2}. Let $M$ be a module as
  specified in \eqclbl{iv}; \thmref{vanishing} yields $\Gdim{M}<
  \infty$. By a result of Holm \thmcite[3.2]{HHl04c}, $R$ is
  Gorenstein if $\id{M}$ is also finite. Now it follows from
  \prpref{fg} that $R$ is Gorenstein.
\end{prf*}

\begin{rmk}
  Let $R$ be a commutative noetherian local ring with residue field
  $k$, and let $E$ denote the injective hull of $k$; it is a
  faithfully injective $R$-module.  A computation based on
  \thmref{stable-Ext-Tor-1} shows that the ranks of $\Stor{i}{k}{k}$
  and $\Sext{i}{k}{k}$ are simultaneously finite:
  \begin{align*}
    \rank_{k}\Stor{i}{k}{k}
    &= \rank_{k}\Stor{i}{k}{\Hom{k}{E}}\\
    &= \rank_{k}\Hom{\Sext{i}{k}{k}}{E}\\
    &= \rank_{k}\Sext{i}{k}{k}\:.
  \end{align*}
  Combined with this equality of ranks, \thmcite[6.4, 6.5, and
  6.7]{LLAOVl07} yield characterizations of regular, complete
  intersection, and Gorenstein local rings in terms of the size of the
  stable homology spaces $\Stor{i}{k}{k}$. For example, $R$ is regular
  if and only if $\Stor{i}{k}{k}=0$ holds for some (equivalently, all)
  $i\in\ZZ$, and $R$ is Gorenstein if and only if $\Stor[R]{i}{k}{k}$
  has finite rank for some (equivalently, all) $i\in\ZZ$.
\end{rmk}

\subsection*{\textbf{\em Commutative rings}} 
If $R$ is commutative and $\p$ is a prime ideal in $R$, then the local
ring $R_\p$ is a flat $R$-algebra, and for an $R_\p$-module $M$ it
follows from \ref{S-structure} that the stable homology modules
$\Stor{i}{M}{N}$ are $R_\p$-modules.

\begin{lem}
  \label{lem:localization}
  Let $R$ be a commutative noetherian ring and let $N$ be an
  $R$-module; let $\p$ be a prime ideal in $R$ and let $M$ be an
  $R_\p$-module. For every $i\in\ZZ$ there is a natural isomorphism of
  $R_\p$-modules,
  \begin{equation*}
    \Stor[R]{i}{M}{N} \dis \Stor[R_\p]{i}{M}{N_\p}\:.
  \end{equation*}
  Hence, $\Stor[R]{i}{-}{N}=0$ implies $\Stor[R_\p]{i}{-}{N_\p}=0$ for
  all prime ideals $\p$ in~$R$.
\end{lem}

\begin{prf*}
  Let $N \qra I$ be an injective resolution. Let $P \qra M$ be a
  projective resolution over $R_\p$; it is a flat resolution of $M$ as
  an $R$-module.  The second isomorphism below follows from
  \prpref{associative}.
  \begin{align*}
    \btp{P}{I}
    &\dis \btp{(\btp[R_\p]{P}{R_\p})}{I}\\
    &\dis \btp[R_\p]{P}{(\btp{R_\p}{I})}\\
    &\dis \btp[R_\p]{P}{(\tp{R_\p}{I})}
  \end{align*}
  The computation gives $\btp{P}{I} \is \btp[R_\p]{P}{I_\p}$;
  similarly one gets $\tp{P}{I} \is \tp[R_\p]{P}{I_\p}$. Now
  \eqref{ts} and the Five Lemma yield $\ttp{P}{I} \is
  \ttp[R_\p]{P}{I_\p}$, and the desired isomorphisms follow as $N_\p
  \to I_\p$ is an injective resolution by Matlis Theory.
\end{prf*}

The proof of the next result is similar. Compared to
\lemref{localization} the noetherian hypothesis on $R$ has been
dropped, as it was only used to invoke Matlis Theory.
\pagebreak
\begin{lem}
  Let $R$ be a commutative ring and let $M$ be an $R$-module; let $\p$
  be a prime ideal in $R$ and let $N$ be an $R_\p$-module. For every
  $i\in\ZZ$ there is a natural isomorphism of $R_\p$-modules,
  \begin{equation*}
    \Stor[R]{i}{M}{N} \dis \Stor[R_\p]{i}{M_\p}{N}\:.
  \end{equation*}
  Hence, $\Stor[R]{i}{M}{-}=0$ implies $\Stor[R_\p]{i}{M_\p}{-}=0$ for
  all prime ideals $\p$~in~$R$.\qed
\end{lem}

\begin{thm}
  \label{thm:finiteid}
  Let $R$ be a commutative noetherian ring and let $N$ be a finitely
  generated $R$-module. If\, $\Stor{i}{-}{N}=0$ holds for some
  $i\in\ZZ$, then $\id[R_\p]{N_\p}$ is finite for every prime ideal
  $\p$ in $R$.
\end{thm}

\begin{prf*}
  From the hypotheses and \lemref{localization} one has
  \smash{$\Stor[R_\p]{i}{-}{N_\p}=0$}. It follows from \prpref{fg}
  that the local ring $R_\p$ is Gorenstein, or $\id[R_\p]N_\p$ is
  finite. However, if $R_\p$ is Gorenstein, then vanishing of
  \smash{$\Stor[R_\p]{i}{-}{N_\p}=0$} implies that $\pd[R_\p]{N_\p}$
  is finite by \corref{Gor3} and \thmref{finitepd}, and then
  $\id[R_\p]{N_\p}$ is finite.
\end{prf*}

The next corollary is now immediate per the remarks in
\ref{ipg:locally}.

\begin{cor}
  \label{cor:cm}
  A commutative noetherian ring $R$ is Cohen--Macaulay if there is a
  finitely generated $R$-module $N\ne 0$ with $\Stor{i}{-}{N}=0$ for
  some $i\in\ZZ$.\qed
\end{cor}

\begin{cor}
  \label{cor:finiteid}
  Let $R$ be a commutative noetherian ring of finite Krull
  dimension. For a finitely generated $R$-module $N$, the following
  conditions are equivalent.
  \begin{eqc}
  \item $\id{N}$ is finite.
  \item $\Stor[R]{i}{-}{N}=0$ for all $i\in\ZZ$.
  \item $\Stor[R]{i}{-}{N}=0$ for some $i\in\ZZ$.
  \end{eqc}
\end{cor}

\begin{prf*}
  The implications
  \eqclbl{i}$\implies$\eqclbl{ii}$\implies$\eqclbl{iii} are clear; see
  \ref{simplifications 2}. Part \eqclbl{i} follows from \eqclbl{iii}
  as $\id{N}$ equals $\sup\setof{\id[R_\p]{N_\p}}{\p\text{ is a prime
      ideal in }R} \le \dimR$.
\end{prf*}

\begin{rmk}
  \label{rmk:id}
  We do not know if the assumption of finite Krull dimension in
  \corref{finiteid} is necessary. By \thmref{finiteid}, vanishing of
  $\Stor{}{-}{N}$ implies that $N$ is locally of finite injective
  dimension, but that does not imply finite injective dimension over
  $R$: Just consider a Gorenstein ring $R$ of infinite Krull
  dimension. On the other hand, we do not know if $\Stor{}{-}{R}$
  vanishes for such a ring.
\end{rmk}

\section{Comparison to Tate homology}
\label{sec:tate}

\noindent
In this section we compare stable homology to Tate homology.  We
parallel some of the findings of Avramov and Veliche
\cite{LLAOVl07}. First we recall a few definitions.

\begin{ipg}
  \label{ipg:gproj}
  An acyclic complex $T$ of projective $\Rop$-modules is called
  \emph{totally acyclic} if $\Hom[\Rop]{T}{P}$ is acyclic for every
  projective $\Rop$-module $P$; cf.~\ref{ipg:tateflat}.  A
  \emph{complete projective resolution} of an $\Rop$-module $M$ is a
  diagram $T \xra{\varpi} P \qra M$, where $T$ is a totally acyclic
  complex of projective $\Rop$-modules, $P\qra M$ is a projective
  resolution, and $\varpi_i$ is an isomorphism for $i \gg 0$; see
  \seccite[2]{OVl06}.

  Let $M$ be an $\Rop$-module with a complete projective resolution $T
  \to P \to M$. For an $R$-module $N$, the \emph{Tate homology} of $M$
  and $N$ over $R$ are the $\kk$-modules $\Ttor{i}{M}{N}=
  \H[i]{\tp{T}{N}}$ for $i\in\ZZ$; see Iacob~\seccite[2]{AIc07}.
\end{ipg}

\begin{ipg}
  \label{ipg:gp}
  An $\Rop$-module has a complete projective resolution if and only if
  it has finite Gorenstein projective dimension; see
  \thmcite[3.4]{OVl06}.

  If $R$ is noetherian and $M$ is a finitely generated $\Rop$-module
  with a complete projective resolution, then $M$ has finite
  G-dimension and it has a complete projective resolution $T \to P
  \to M$ with $T$ and $P$ degree-wise finitely generated and $T\to P$
  surjective; see \cite[2.4.1]{OVl06} and \thmcite[3.1]{LLAAMr02}.
\end{ipg}

\begin{lem}
  \label{lem:Gproj-Gflat}
  Let $M$ be an $\Rop$-module that has a complete projective
  resolution $T \to P \to M$. The following conditions are equivalent.
  \begin{eqc}
  \item $T\otimes_{R}E$ is acyclic for every injective $R$-module $E$.
  \item There are isomorphisms of functors $\Stor{i}{M}{-} \is
    \Ttor{i}{M}{-}$ for all $i\in\ZZ$.
  \end{eqc}
\end{lem}

\begin{prf*}
  Assume that $\tp{T}{E}$ is acyclic for every injective $R$-module
  $E$. The pair $(T,P)$ is then a Tate flat resolution of $M$, see
  \ref{ipg:tateflat}, and it follows from \thmref{Stor-Tateflat} that
  the functors $\Stor{i}{M}{-}$ and $\Ttor{i}{M}{-}$ are isomorphic
  for all $i\in\ZZ$.  For the converse, let $E$ be an injective
  $R$-module. By \ref{simplifications 2} one then has
  \begin{equation*}
    0 = \Stor{i}{M}{E} \is \Ttor{i}{M}{E}
  \end{equation*}
  for all $i\in\ZZ$, and hence $\H{\tp{T}{E}}=0$.
\end{prf*}

\begin{thm}
  \label{thm:Gproj-Gflat-fg}
  Let $R$ be noetherian, and let $M$ be a finitely generated
  $\Rop$-module that has a complete projective resolution. There are
  isomorphisms of functors
  \begin{equation*}
    \Stor{i}{M}{-} \dis \Ttor{i}{M}{-}\quad\text{for all } i\in\ZZ\:.
  \end{equation*}
\end{thm}

\begin{prf*}
  The module $M$ has a complete projective resolution $T \to P \to M$
  with $T$ and $P$ degree-wise finitely generated; see
  \pgref{ipg:gp}. The isomorphisms \eqref{gflatTate} show that the
  complex $\tp{T}{E}$ is acyclic for every injective $R$-module $E$,
  and \lemref{Gproj-Gflat} finishes the argument.
\end{prf*}

\begin{rmk}
  The isomorphisms of homology modules in \thmref{Gproj-Gflat-fg}
  actually follow from \textsl{one} isomorphism in $\D[\kk]$, but this
  is unapparent in the proof, which rests on \thmref{Stor-Tateflat}.
  The finitely generated module $M$ has a complete projective
  resolution $T \to L \to M$ with $T$ and $L$ degree-wise finitely
  generated and \mbox{$T \to L$} surjective; see \pgref{ipg:gp}.  Thus the
  kernel $K = \Ker(T \to L)$ is a bounded above complex of projective
  modules. Given a module $N$, let $C$ be the cone as one would construct in
  \ref{ipg:C}. By \prpref{btp}(a) the complex $\btp{K}{C}$ is acyclic,
  so that \lemref{tTor} and \ref{btp}(c) yield $\susp^{-1}\ttpp{L}{I}
  \qis \susp^{-1}\btpp{L}{C} \qis \susp^{-1}\btpp{T}{C}$.  As in the
  proof above, $\tp{T}{E}$ is acyclic, so \lemref{N} gives
  $\susp^{-1}\btpp{T}{C} \qis \tp{T}{N}$, and combining these
  isomorphisms in $\D[\kk]$ gives the desired one.
\end{rmk}

Without extra assumptions on the ring, we do not know if stable
homology agrees with Tate homology whenever the latter is defined. In
general, the relation between stable homology and Tate homology is
tied to an unresolved problem in Gorenstein homological
algebra. \thmref{Gproj-Gflat} explains how.

\begin{ipg}
  \label{ipg:gpgf}
  An $\Rop$-module $G$ is called \emph{Gorenstein projective} if there
  exists a totally acyclic complex $T$ of projective $\Rop$-modules
  with $\Co[0]{T} \is G$; see \ref{ipg:gproj}. Similarly, an
  $\Rop$-module $G$ is called \emph{Gorenstein flat} if there exists a
  totally acyclic complex $T$ of flat $\Rop$-modules with $\Co[0]{T}
  \is G$; see \ref{ipg:tateflat}.
\end{ipg}
\pagebreak
\begin{thm}
  \label{thm:Gproj-Gflat}
  The following conditions on $R$ are equivalent.
  \begin{eqc}
  \item Every Gorenstein projective $\Rop$-module is Gorenstein flat.
  \item For every $\Rop$-module $M$ that has a complete projective
    resolution there are isomorphisms of functors $\Stor{i}{M}{-} \is
    \Ttor{i}{M}{-}$ for all $i\in\ZZ$.
  \end{eqc}
\end{thm}

\begin{prf*}
  Assume that every Gorenstein projective $\Rop$-module is Gorenstein
  flat. Let $T \to P \to M$ be a complete projective resolution. It
  follows that $T$ is a totally acyclic complex of flat modules; see
  Emmanouil~\thmcite[2.2]{IEm12}. Thus stable homology and Tate
  homology coincide by \lemref{Gproj-Gflat}.

  For the converse, let $M$ be a Gorenstein projective $\Rop$-module
  and let $T$ be a totally acyclic complex of projective
  $\Rop$-modules with $M \is \Co[0]{T}$.  Since there are isomorphisms
  of functors $\Stor{i}{M}{-} \is \Ttor{i}{M}{-}$ for all $i\in\ZZ$,
  it follows from \ref{simplifications 2} that $\tp{T}{E}$ is acyclic
  for every injective $R$-module $E$. Thus $T$ is a totally acyclic
  complex of flat $\Rop$-modules, and so $M$ is Gorenstein flat.
\end{prf*}

\begin{rmk}
  \label{rmk:Gproj-Gflat}
  As Holm notes \prpcite[3.4]{HHl04a}, the obvious way to achieve that
  every Gorenstein projective $\Rop$-module is Gorenstein flat is to
  ensure that (1) the Pontryagin dual of every injective $R$-module is
  flat, and (2) that every flat $\Rop$-module has finite projective
  dimension. The first condition is satisfied if $R$ is left coherent,
  and the second is discussed in \rmkref{pdfm}.  A description of the
  rings over which Gorenstein projective modules are Gorenstein flat
  seems elusive; see \seccite[2]{IEm12}.
\end{rmk}

\subsection*{\textbf{\em Complete homology}} 
In his thesis, Triulzi considers the \emph{J-completion} of the
homological functor $\Tor{}{M}{-} =
\setof{\Tor{i}{M}{-}}{i\in\ZZ}$. His construction is similar to
Mislin's P-completion of covariant Ext and Nucinkis' I-completion of
contravariant Ext. The resulting homology theory is called
\emph{complete homology.} Like stable homology, it is a generalization
of Tate homology. We compare these two generalizations in
\cite{CCLPa}. From the point of view of stable homology, it is
interesting to know when it agrees with complete homology, because the
latter has a universal property. In this direction the main results in
\cite{CCLPa} are that these two homology theories agree over
Iwanaga-Gorenstein rings, and for finitely generated modules over Artin
algebras and complete commutative local rings. Moreover, the two
theories agree with Tate homology, whenever it is defined, under the
exact same condition; that is, if and only if every Gorenstein
projective module is Gorenstein flat; see \thmref{Gproj-Gflat}.

\appendix
\section*{Appendix}
\label{stable cohomology}
\stepcounter{section}
\noindent
We start by recalling the definition of stable cohomology.

\begin{ipg}
  \label{ipg:bHom}
  Let $X$ and $Y$ be $R$-complexes, following \cite{LLAOVl07,FGc92} we
  denote by $\bHom{X}{Y}$ the subcomplex of $\Hom{X}{Y}$ with degree
  $n$ term
  \begin{equation*}
    \bHom{X}{Y}_{n} = \coprod_{i\in\ZZ}\Hom{X_{i}}{Y_{i+n}}\:.
  \end{equation*}
  It is called the \emph{bounded} Hom complex, and the quotient
  complex
  \begin{equation*}
    \SHom{X}{Y} = \Hom{X}{Y}/\bHom{X}{Y}
  \end{equation*}
  is called the \emph{stable} Hom complex.

  For $R$-modules $M$ and $N$ with projective resolutions $P_M \qra M$
  and $P_N \qra N$, the $\kk$-modules
  \begin{equation*}
    \bExt{i}{M}{N} = \H[-i]{\bHom{P_{M}}{P_{N}}}\:,
  \end{equation*}
  are called the \emph{bounded cohomology} of $M$ and $N$ over $R$,
  and the \emph{stable cohomology} modules of $M$ and $N$ over $R$ are
  \begin{equation*}
    \Sext{i}{M}{N} = \H[-i]{\SHom{P_{M}}{P_{N}}}\:.
  \end{equation*}

  Avramov and Veliche~\cite{LLAOVl07} use the notation
  $\Text{i}{M}{N}$ for the stable cohomology; this notation is
  standard for Tate cohomology, which coincides with stable cohomology
  whenever the former is defined; see \corcite[2.4]{LLAOVl07}.
\end{ipg}

\begin{prp}
  \label{prp:bhom}
  Let $X$ and $Y$ be $R$-complexes.
  \begin{prt}
  \item I\fsp[.4] $X$ or $Y$ is bounded above, and $\Hom{X_i}{Y}$ is
    acyclic for all $i$, then the complex $\bHom{X}{Y}$ is acyclic.
  \item I\fsp[.4] $X$ or $Y$ is bounded below, and $\Hom{X}{Y_i}$ is
    acyclic for all $i$, then the complex $\bHom{X}{Y}$ is acyclic.
  \end{prt}
\end{prp}

\begin{prf*}
  Similar to the proof of \prpref{btp}.
\end{prf*}

\subsection*{\textbf{\em Standard isomorphisms}} We study composites
of the functors $\widebar{\rm Hom}$ and $\widebar{\otimes}$. To the
extent possible, we establish analogs of the standard isomorphisms for
composites of Hom and $\otimes$; see \seccite[II.5 and
VI.5]{careil}. There seems to be no analog of Hom-tensor adjunction
\prpcite[II.5.2]{careil}.

The setup is the same for Propositions
\ref{prp:associative}--\ref{prp:bhev}; namely:

\begin{ipg}
  \label{appendix notation}
  Let $X$ be a complex of $\Rop$-modules, let $Y$ be a complex of
  $(R,\Sop)$-bimodules, and let $Z$ be a complex of $S$-modules.
\end{ipg}

Under finiteness conditions, the unbounded tensor product is
associative.

\begin{prp}
  \label{prp:associative} For complexes as in \ref{appendix notation},
  under either of the following conditions%
  \begin{itemlist}
  \item $X$ and $Z$ are complexes of finitely presented modules
  \item $Y$ is a bounded complex
  \end{itemlist}
  there is an isomorphism o\fsp\ $\kk$-complexes,
  \begin{equation*}
    \btp[S]{(\btp{X}{Y})}{Z} \;\lra\; \btp{X}{(\btp[S]{Y}{Z})}\:,
  \end{equation*}
  and it is functorial in $X$, $Y$, and $Z$.
\end{prp}

\begin{prf*}
  For every $n\in\ZZ$ one has,
  \begin{align*}
    \btpp[S]{\btpp{X}{Y}}{Z}_n &= \prod_{i\in\ZZ}
    \tp[S]{\left(\textstyle\prod_{j\in\ZZ}\tp{X_j}{Y_{i-j}}\right)}{Z_{n-i}}\\[-.6\baselineskip]
    \intertext{and}\\[-1.65\baselineskip]
    \btpp{X}{\btpp[S]{Y}{Z}}_n 
    &= \prod_{j\in\ZZ}
    \tp{X_j}{\left(\textstyle\prod_{i\in\ZZ}\tp[S]{Y_{i-j}}{Z_{n-i}}\right)}\:,
  \end{align*}
  and from each of these modules there is a canonical homomorphism to
  \begin{equation*}
    \prod_{i\in\ZZ}\prod_{j\in\ZZ} \tp[S]{\tpp{X_j}{Y_{i-j}}}{Z_{n-i}}
    \dis \prod_{j\in\ZZ}\prod_{i\in\ZZ} \tp{X_j}{\tpp[S]{Y_{i-j}}{Z_{n-i}}}\:.
  \end{equation*}
  If $Y$ is bounded, then these homomorphisms are isomorphisms as the
  inner most products are finite.  Recall, e.g.\ from
  \thmcite[3.2.22]{rha}, that the functor $\tp{M}{-}$ commutes with
  products if $M$ is finitely presented.  Thus, the canonical
  homomorphisms are isomorphisms when $X$ and $Z$ are complexes of
  finitely presented modules. It is straightforward to verify that
  these isomorphisms commute with the differentials and form an isomorphism of
  complexes.
\end{prf*}

The model for the following \emph{swap} isomorphism is \cite[ex.\
II.4]{careil}. The proof is similar to the proof of
\prpref{associative} and uses that the functor $\Hom{M}{-}$ commutes
with coproducts if $M$ is a finitely generated $R$-module.

\begin{prp}
  \label{prp:bswp}
  For complexes as in \ref{appendix notation}, under either of the
  following conditions%
  \begin{itemlist}
  \item $X$ and $Z$ are complexes of finitely generated modules
  \item $Y$ is a bounded complex
  \end{itemlist}
  there is an isomorphism o\fsp\ $\kk$-complexes,
  \begin{equation*}
    \bHom[\Rop]{X}{\bHom[S]{Z}{Y}} \;\lra\; \bHom[S]{Z}{\bHom[\Rop]{X}{Y}}\:,
  \end{equation*}
  and it is functorial in $X$, $Y$, and $Z$.\qed
\end{prp}

The next result is an analog of \prpcite[VI.5.2 and VI.5.3]{careil};
it is used below to establish a duality between
stable homology and stable cohomology.

\begin{prp}
  \label{prp:bhev}
  Let $X$, $Y$, and $Z$ be as in \ref{appendix notation}, and assume
  further that $X$ is a complex of finitely presented
  $\Rop$-modules. There is a morphism o\fsp\ $\kk$-complexes,
  \begin{equation*}
    \btp{X}{\Hom[S]{Y}{Z}} \;\lra\; \Hom[S]{\bHom[\Rop]{X}{Y}}{Z}\:,
  \end{equation*}
  and it is functorial in $X$, $Y$, and $Z$.  Furthermore, it is an
  isomorphism if $X$ is a complex of projective modules or $Z$ is a
  complex of injective modules.
\end{prp}

\begin{prf*}
  For every $n\in\ZZ$ one can compute as follows,
  \begin{align*}
    \btpp{X}{\Hom[S]{Y}{Z}}_{n}
    &= \prod_{i\in\ZZ}\tpp{X_{i}}{\Hom[S]{Y}{Z}_{n-i}}\\
    &= \prod_{i\in\ZZ} \tpp{X_{i}}{\dprod_{j\in\ZZ}\Hom[S]{Y_{j}}{Z_{n-i+j}}}\\
    &\is
    \prod_{i\in\ZZ}\prod_{j\in\ZZ}\tpp{X_{i}}{\Hom[S]{Y_{j}}{Z_{n-i+j}}}\:,
  \end{align*}
  where the isomorphism holds as the module $X_{i}$ is finitely
  presented for every $i\in\ZZ$; see \thmcite[3.2.22]{rha}. On the
  other hand, for every $n\in\ZZ$ one has
  \begin{align*}
    \Hom[S]{\bHom[\Rop]{X}{Y}}{Z}_{n}
    &= \prod_{h\in\ZZ} \Hom[S]{\bHom[\Rop]{X}{Y}_{h}}{Z_{n+h}}\\
    &=\prod_{h\in\ZZ} \Hom[S]{\dcoprod_{i\in\ZZ}\Hom[\Rop]{X_{i}}{Y_{i+h}}}{Z_{n+h}}\\
    &\is \prod_{h\in\ZZ}\prod_{i\in\ZZ} \Hom[S]{\Hom[\Rop]{X_{i}}{Y_{i+h}}}{Z_{n+h}}\\
    &= \prod_{i\in\ZZ}\prod_{j\in\ZZ}
    \Hom[S]{\Hom[\Rop]{X_{i}}{Y_{j}}}{Z_{n-i+j}}\:.
  \end{align*}
  Now set $(\theta_{XYZ})_{n} =
  \prod_{i\in\ZZ}\prod_{j\in\ZZ}(-1)^{i(n-j)}\theta_{X_{i}Y_{j}Z_{n-i+j}}$,
  where
  \begin{equation*}
    \dmapdef{\theta_{X_{i}Y_{j}Z_{n-i+j}}}{\tp{X_{i}}{\Hom[S]{Y_{j}}{Z_{n-i+j}}}}{
      \Hom[S]{\Hom[\Rop]{X_{i}}{Y_{j}}}{Z_{n-i+j}}}
  \end{equation*}
  is the homomorphism of $\kk$-modules given by
  $\theta_{X_{i}Y_{j}Z_{n-i+j}}(x\otimes\psi)(\phi) = \psi\phi(x)$.
  It is straightforward to verify that $\theta_{XYZ}$ is a morphism of
  $\kk$-complexes and functorial in $X$, $Y$, and $Z$.

  Finally, if each module $X_i$ is projective or each module
  $Z_{n-i+j}$ is injective, then it follows from \prpcite[VI.5.2 and
  VI.5.3]{careil} that $\theta_{X_{i}Y_{j}Z_{n-i+j}}$ is invertible
  for all $i,j,n\in\ZZ$, and so the morphism $\theta_{XYZ}$ is
  invertible.
\end{prf*}

\begin{lem}
  \label{lem:stable-Ext-Tor-1}
  Let $P$ be a bounded below complex of finitely generated projective
  $\Rop$-modules and let $X$ be a complex of $\Rop$-modules with $\H{X}$
  bounded. For every injective $\kk$-module $E$, there is an
  isomorphism in the derived category $\D[\kk]$:
  \begin{equation*}
    \ttp{P}{\Hom[\kk]{X}{E}} \lra \susp\Hom[\kk]{\SHom[\Rop]{P}{X}}{E}\;,
  \end{equation*}
  and it is functorial in $P$, $X$, and $E$.
  \end{lem}

\begin{prf*}
  Set $(-)^{\vee}=\Hom[\kk]{-}E$.  In the commutative square of $\kk$-complexes
  \begin{equation*}
    \xymatrix{
      \tp{P}{X^{\vee}} \ar[r]^-{\tau} \ar[d]^-\theta
      & \btp{P}{X^{\vee}} \ar[d]_-\is^-\a\\
      \Hom[\Rop]{P}{X}^{\vee} \ar[r]^-{\vartheta^\vee} & \bHom[\Rop]{P}{X}^{\vee}
    }
  \end{equation*}
  each horizontal morphism is (the dual of) a canonical embedding. The
  vertical map $\a$ is the isomorphism from \prpref{bhev}.  The
  morphism $\theta$ is the standard evaluation map; it is a
  quasi-isomorphism by~\cite[4.4(I)]{LLAHBF91}. The square induces a
  morphism of triangles in the homotopy category:
  \begin{equation*}
    \xymatrix{
      \tp{P}{X^{\vee}} \ar[r]^-{\tau} \ar[d]^-\theta_-{\qis}
      & \btp{P}{X^{\vee}} \ar[d]_-\is^-\a \ar[r] & \Cone{\tau}
      \ar[d]^-{\gamma}
      \ar[r] & \susp\tpp{P}{X^{\vee}} \ar[d]^-{\susp\theta}_-\qis\\
      \Hom[\Rop]{P}{X}^{\vee} \ar[r]^-{\vartheta^\vee} &
      \bHom[\Rop]{P}{X}^{\vee} \ar[r]
      & \Cone{\vartheta^\vee} \ar[r] &
      \susp\Hom[\Rop]{P}{X}^{\vee}}
  \end{equation*}
  The construction of $\gamma$ is functorial in all three arguments,
  and it is a quasi-isomorphism because $\a$ and $\theta$ are
  quasi-isomorphisms. Recall, say from \cite[III.3.4--5]{GelMan}, that
  there are natural quasi-isomorphisms
  \begin{gather*}
    \Cone{\tau} \dqis \Coker{\tau} \deq \ttp{P}{X^{\vee}}\qand\\
    \Cone{\vartheta^\vee} \dis \susp(\Cone{\vartheta})^\vee
    \dqis \susp(\Coker{\vartheta})^\vee \deq \susp\SHom[\Rop]{P}{X}^\vee\;.
  \end{gather*}
  They yield the desired isomorphism in the derived category
  \begin{equation*}
    \ttp{P}{\Hom[\kk]{X}{E}} \lra \susp\Hom[\kk]{\SHom[\Rop]{P}{X}}{E}\;.
  \end{equation*}
  With regard to functoriality of $\gamma$, notice that a morphism between
  arguments, $P \to P_1$ say, induces the solid commutative square in
  the following diagram.
  \begin{equation*}
    \xymatrix@C=.5pc@R=1.8pc{
      \Coker{\tau} \ar@{-->}[rrrr] \ar@{-->}[ddd]
      & & & & \Coker{\tau_1} \ar@{-->}[ddd]\\
      & \Cone{\tau} \ar[rr] \ar[d]^-\gamma \ar[lu]^-\qis
      && \Cone{\tau_1} \ar[ru]_\qis \ar[d]^-{\gamma_1}  \\
      & \Cone{\vartheta^\vee} \ar[rr] \ar[ld]_-\qis
      && \Cone{(\vartheta_1^\vee)} \ar[rd]^\qis\\
      \susp(\Coker{\vartheta)^\vee} \ar@{-->}[rrrr]
      & & & & \susp(\Coker{\vartheta_1})^\vee
      }
  \end{equation*}
  Commutativity in the derived category of the dashed square is now a
  consequence. Functoriality in the other arguments is handles similarly.
\end{prf*}

The Lemma immediately yields a useful duality.

\begin{thm}
  \label{thm:stable-Ext-Tor-1}
  Let $M$ and $N$ be $\Rop$-modules and assume that $M$ has a
  degree-wise finitely generated projective resolution. For every
  injective $\kk$-module $E$ and for every $i\in\ZZ$ there is an
  isomorphism of\, $\kk$-modules,
  \begin{equation*}
    \Hom[\kk]{\Sext[\Rop]{i}{M}{N}}{E} \dis \Stor{i}{M}{\Hom[\kk]{N}{E}}\:,
  \end{equation*}
  and it is functorial in $M$, $N$, and $E$. \qed
\end{thm}

To prove the next two results one proceeds as in the proof of
\prpref{bhev}.

\begin{prp}
  \label{prp:tev}
  Let $X$ be a complex of finitely generated $R$-modules and let $Y$
  and $Z$ be as in \ref{appendix notation}.  There is a morphism
  o\fsp\ $\kk$-complexes
  \begin{equation*}
    \tp[S]{\bHom{X}{Y}}{Z} \;\lra\; \bHom{X}{\tp[S]{Y}{Z}}\:,
  \end{equation*}
  and it is functorial in $X$, $Y$, and $Z$. Furthermore, it is an
  isomorphism under each of the following conditions
  \begin{itemlist}
  \item $Z$ is a complex of finitely generated projective modules
  \item $X$ is a complex of finitely presented modules and $Z$ is a
    complex of flat~modules
  \item $X$ is a complex of projective modules \qed
  \end{itemlist}
\end{prp}

\begin{prp}
  \label{prp:tevb}
  Let $X$ be a complex of $R$-modules, let $Y$ be a complex of
  $(R,\Sop)$-bimodules, and let $Z$ be a complex of finitely presented
  $S$-modules. There is a morphism o\fsp\ $\kk$-complexes,
  \begin{equation*}
    \btp[S]{\Hom{X}{Y}}{Z} \;\lra\; \Hom{X}{\btp[S]{Y}{Z}}\:,
  \end{equation*}
  and it is functorial in $X$, $Y$, and $Z$. Furthermore, it is an
  isomorphism if $X$ or $Z$ is a complex of projective modules.\qed
\end{prp}

\subsection*{\textbf{\em Pinched tensor products}}
Christensen and Jorgensen devised in \cite{LWCDAJ14} a \emph{pinched}
tensor product, $\ptp{-}{-}$, to compute Tate homology. In view of
\thmref{Stor-Tateflat} their proof of \thmcite[3.5]{LWCDAJ14} applies
verbatim to yield the next result; we refer the reader to
\cite{LWCDAJ14} for the definition of the pinched tensor product.

\begin{thm}
  \protect\pushQED{\qed}%
  \label{thm:Stor-pinch}
  Let $M$ be an $\Rop$-module that has a Tate flat resolution $(T,F)$,
  let $A$ be an acyclic complex of $R$-modules and set
  $N=\Co[0]{A}$. For every $i\in\ZZ$, there is an isomorphism o\fsp\
  $\kk$-modules
  \begin{equation*}
    \Stor{i}{M}{N} \dis \H[i]{\ptp{T}{A}}\:. \qedhere
  \end{equation*}
\end{thm}

The next corollary is an analogue of \cite[Corollary 4.10]{LWCDAJ14}.

\begin{cor}
  \label{cor:Stor-pinch}
  Let $R$ be commutative and let $M$ and $N$ be Gorenstein flat
  $R$-modules with corresponding totally acyclic complexes of flat
  modules $T$ and $T'$, respectively. If $\Stor{i}{M}{N}=0$ holds for
  all $i\in\ZZ$, then $\ptp{T}{T'}$ is an acyclic complex of flat
  $R$-modules, and the following statements are equivalent:
  \begin{eqc}
  \item $\ptp{T}{T'}$ is a totally acyclic complex of flat
    $R$-modules.
  \item $\Stor{i}{M}{\tp{N}{E}}=0$ holds for every injective
    $R$-module $E$ and all $i\in\ZZ$.
  \end{eqc}
  When these conditions hold, $M\otimes_{R}N$ is a Gorenstein flat
  $R$-module and $\ptp{T}{T'}$ is a corresponding totally acyclic
  complex of flat $R$-modules.
\end{cor}

\begin{prf*}
  It follows from the definition of pinched tensor products that
  $\ptp{T}{T'}$ is a complex of flat $R$-modules, and if
  $\Stor{i}{M}{N}=0$ holds for all $i\in\ZZ$, then the complex is
  acyclic by \thmref{Stor-pinch}. It is totally acyclic if and only if
  $\tp{\ptpp{T}{T'}}{E} \is \ptp{T}{\tpp{T}{E}}$ is acyclic for every
  injective $R$-module $E$; that is, if and only if
  $\Stor{i}{M}{\tp{N}{E}}=0$ holds for every injective $R$-module $E$
  and all $i\in\ZZ$. Finally, it follows from the definition of
  pinched tensor products that there is an isomorphism $\tp{M}{N} \is
  \Co[0]{\ptp{T}{T'}}$.
\end{prf*}


\providecommand{\arxiv}[2][AC]{\mbox{\href{http://arxiv.org/abs/#2}{\sf arXiv:#2 [math.#1]}}}
\providecommand{\MR}[1]{\mbox{\href{http://www.ams.org/mathscinet-getitem?mr=#1}{#1}}}
\renewcommand{\MR}[1]{\mbox{\href{http://www.ams.org/mathscinet-getitem?mr=#1}{#1}}}
\providecommand{\bysame}{\leavevmode\hbox to3em{\hrulefill}\thinspace}
\providecommand{\MR}{\relax\ifhmode\unskip\space\fi MR }
\providecommand{\MRhref}[2]{%
  \href{http://www.ams.org/mathscinet-getitem?mr=#1}{#2}
}
\providecommand{\href}[2]{#2}

\end{document}